\documentclass[reqno,11pt]{amsart}
\usepackage[top=1.2in,bottom=1.0in,left=1.0in,right=1.0in]{geometry}
\usepackage{amsmath}
\usepackage{amssymb}
\usepackage{mathrsfs}
\usepackage{times}
\usepackage{bm}
\usepackage{latexsym}
\usepackage{graphicx}
\usepackage{indentfirst}
\usepackage{graphicx}
\usepackage{epsfig}
\usepackage{psfrag}
\usepackage{subfigure}
\usepackage{url}
\usepackage{caption}
\usepackage{float}
\usepackage{stfloats}
\usepackage{fancyhdr}
\usepackage{dsfont}
\usepackage{amsfonts}
\usepackage{enumerate}
\usepackage{epstopdf}
\usepackage{color}
\usepackage{amssymb,amsmath,amsthm,mathrsfs}
\usepackage{booktabs}
\usepackage{threeparttable}
\usepackage[titletoc]{appendix}
\usepackage{algorithm}
\usepackage{algorithmic}
\usepackage{multirow}
\usepackage{array}
\usepackage{diagbox}
\usepackage[misc]{ifsym}
\usepackage{boxedminipage}
\usepackage{enumerate}

    \numberwithin{equation}{section}%
    \numberwithin{table}{section}%
    \numberwithin{figure}{section}
\allowdisplaybreaks

\newtheorem{lemma}{Lemma}[section]
\newtheorem{remark}{Remark}[section]

\newtheorem{theorem}{Theorem}[section]

\DeclareMathAlphabet{\mathmatrix}{OT1}{ptm}{b}{n}
\DeclareMathAlphabet{\mathvector}{OT1}{ptm}{bx}{it}

\def\d{\mathrm{d}}
\def\RR{\mathbb{R}}

\def\det{\mathrm{det}}

\def \bx{\bm x}

\begin{document}

\title[]{ Theoretical  analysis   of  the extended cyclic reduction algorithm  }

\author{Xuhao Diao}
\address[Xuhao Diao]{ School of Mathematical Sciences, Peking University, Beijing 100871, China.}
\email{diaoxuhao@pku.edu.cn}

\author{Jun Hu}
\address[Jun Hu]{ School of Mathematical Sciences, Peking University, Beijing 100871, China.} \email{hujun@math.pku.edu.cn}

\author{Suna Ma}
\address[Suna Ma \Letter]{ Nanjing University of Posts and Telecommunications, Nanjing 210023, China.} \email{masuna@njupt.edu.cn}

\keywords{ extended cyclic reduction algorithm, forward error analysis, block-tridiagonal linear systems,  principal submatrix, eigenpolynomial}

\subjclass[2000]{Primary 65N35, 65F05, 65N12}

\maketitle

\begin{abstract}

The extended cyclic reduction algorithm developed by Swarztrauber in 1974 was used to solve the
block-tridiagonal linear system.
The paper fills in the gap of theoretical  results concerning the zeros of matrix polynomial  $B_{i}^{(r)}$ with respect to a tridiagonal matrix which are computed by Newton's method in  the extended cyclic reduction algorithm.
Meanwhile, the forward error analysis of the extended cyclic reduction algorithm for solving the block-tridiagonal system is studied. To achieve the two aims, the critical point is to find out that the zeros of  matrix polynomial $B_{i}^{(r)}$  are eigenvalues of a  principal submatrix of the coefficient matrix. 

\end{abstract}

\section{introduction}

Tridiagonal and block-tridiagonal systems  play a fundamental role in matrix computation related to scientific and engineering problems, which  particularly occur in approximation of the
finite difference method for the Poisson equation. Among many algorithms to solve such systems,  there are two basic techniques of the direct methods  which are computationally very fast and require a minimum of storage:
\begin{itemize}
\item Fast Fourier transform, which relies on the knowledge of a certain set of
trigonometric eigenvectors.
\item Cyclic reduction, which relies on the simple block tridiagonal structure
of the coefficient matrix.
\end{itemize}
Cyclic reduction  is an algorithm invented by  Golub and
 Hockney in the mid 1960s for solving linear systems resulting from the finite
difference method for the Poisson equation over a rectangle.
Since then it  received
much attention for its very nice computational features and had a great development \cite{Bank1975,BuzbeeGolub1970,BuzbeeDorr1971,BuzbeeDorr1974,Swarztrauber1974JCP,BankRose1975,
Swarztrauber1973,BiniMeini2009}. Among the algorithms related by Golub, it is one of the most versatile and powerful ever created \cite{BiniMeini2009}.
Afterwards Swarztrauber  extended the cyclic reduction method to linear systems related to the discretization of separable elliptic equations with Dirichlet, Neumann, or periodic boundary conditions \cite{Swarztrauber1974}, which is the so-called  extended cyclic reduction (ECR) algorithm. The ECR algorithm  has been adopted by FISHPACK  which is an efficient FORTRAN subprograms for the solution of separable elliptic partial differential equations  by Adams, Swarztrauber and  Sweet \cite{Swarztrauber1975}.
The ECR algorithm for the discrete system resulting from separable elliptic equations was noted in classic book \cite{saad2003} by  Saad.

Of particular interests are direct methods for  linear systems with the block-tridiagonal matrices resulting from Legendre-Galerkin spectral methods for the constant-coefficient elliptic equations. A  direct approach by the matrix diagonalization method \cite{HaidvogelZang1979}
which was based on the spectral  decomposition of  matrices was presented  for the Legendre-Galerkin approximation of the two and three dimensional Helmholtz equations by Shen in \cite{Shen1994Legendre}, whose complexity is of $\mathcal{O}(N^{d+1})$, where $d=2,3$ and $N$ is the cutoff number of the polynomial expansion in each direction.
 A fast direct two-dimensional Poisson solver, the complexity of which is better than that of the algorithm based on the  matrix diagonalization method,
 was constructed by
 further exploring the matrix structures of the Legendre-Galerkin spectral discretization combinated with the ECR algorithm \cite{Shen1995},  whose complexity is of $\mathcal{O}(N^{2}\log_2 N)$.
Recently, it was extended to fast solve  the three-dimensional Poisson equation in \cite{Legendre1}, whose complexity is of $\mathcal{O}(N^{3}(\log_2 N)^2)$. As mentioned above,
 direct methods based on cyclic reduction yield a quasi-optimal complexity for the systems related to the Legendre-Galerkin spectral discretization. To the best of our knowledge, no related theoretical analysis has been given about the ECR algorithm and this paper is to focus on it.

Cyclic reduction for the block-tridiagonal system  described in Section 2  proceeds by first eliminating half of the variables simultaneously, then half of the remainder, and so on. This process is continued until a system with a single unknown vector  is obtained.  In the implementation, each step generates a block-tridiagonal system with the matrices $I, 2I-(A^{(r)})^2, I$, where $I$ is the identity matrix and the matrix $A^{(r)}= P_{2^r}(A)$  is a polynomial of degree
not greater than $2^r$ with respect to $A$. 
Actually, the matrix  in $r$-th reduction step is
$$A^{(r+1)}=2I-(A^{(r)})^2,$$
which  indicates  the polynomial $P_{2^r}(x)$ satisfies the following recurrence relation
\begin{align*}
  P_1(x)=x, \qquad  P_{2^{r+1}}(x)=2-P^2_{2^{r}}(x).
\end{align*}
This together with the property of the Chebyshev polynomials $T_n(t)$ leads to
\begin{align*}
  P_{2^r}(x)=-2T_{2^r}(-\frac{x}{2}), \qquad r\geq 0.
\end{align*}
Since the zeros of Chebyshev polynomials are available explicitly, the  matrix  polynomial $A^{(r)}$ of matrix $A$ can be expressed in a factorized form  directly  as is shown in \eqref{CRAr}. 
 By contrast, in the implementation of the ECR algorithm, each step generates a block-tridiagonal system with  matrices $A_i^{(r)},B_i^{(r)}, C_i^{(r)}$ (see \eqref{kuaisanjiao} below) which are  polynomials of the  tridiagonal matrix $B$.
 Since matrices $A_i^{(r)}, B_i^{(r)}, C_i^{(r)}$ fill rapidly  as $r$ increases,  and consequently storage requirements of computation  become excessive.
 Instead of storing these matrices, it is to  store the zeros of the corresponding polynomials that represent them. But the explicit expressions of these zeros  are not available like the cyclic reduction \cite{Swarztrauber1974}.

 It is observed from numerical experiments that  Newton's method  is globally convergent only if the zeros of  matrix polynomial $B_i^{(r)}$  are real and simple.
 However, a theoretical proof is missing in literatures. Thus the first  aim of this article is to prove the zeros of  matrix polynomial $B_i^{(r)}$ are real and simple so that the algorithm is  stable and globally convergent in theory. Here and  throughout this paper, the zeros of a matrix polynomial with respect to some matrices mean the zeros of the corresponding scalar valued  polynomial that represents it.
  The second aim  is to give the forward error analysis of  the  ECR algorithm for the block-tridiagonal system.
For this purpose, the block-tridiagonal system is rewritten in the  tensor product formulation
$\big( B\otimes I + I\otimes \mathcal{R}_n \big)X=Y$, which is presented in Section 3. And
   two basic assumptions are made:
(1) Both $B$ and $\mathcal{R}_n$ are symmetric positive definite  tridiagonal matrices, and all subdiagonal entries of $\mathcal{R}_n$ are non-zero.
 (2) The eigenvalues of  $B$ and $\mathcal{R}_n$ satisfy
 \begin{align*}
   \lambda_{\text{max}}(\mathcal{R}_n)\le
    \max_{1\le i\le n}|b_i\pm (|a_i|+|a_{i+1}|) |&\le 1,  \qquad  \lambda_{\text{min}}(\mathcal{R}_n)\geq \mathbf{u}^{1-\epsilon},\\
    \lambda_{\text{max}}(B)&\le 1,  \qquad  \lambda_{\text{min}}(B)\geq \mathbf{u}^{1-\epsilon},
 \end{align*}
 where $ 0<\epsilon<1,$  $\mathbf{u}$ is the unit roundoff.
   Under the two assumptions, the forward error analysis of the ECR algorithm
   in the reduction phase and the back-substitution phase
   is given.

The main contributions of this paper  are of three folds:
\begin{enumerate}[(i)]
\item It is proved that the zeros of  matrix polynomial $B_i^{(r)}$  are real and simple which fills in the gap of theoretical results.
\item It is found out that the matrix polynomial $B_i^{(r)}$  corresponds to the eigenpolynomial of a principal submatrix of $-\mathcal{R}_n$ which is the  critical point of the theoretical analysis throughout the paper. Thanks to this conclusion, the method of bisection (MOB) for the eigenvalues of a symmetric tridiagonal matrix can be applied  to compute the  zeros of  $B_{i}^{(r)}$,  the computational cost and accuracy of which are  quantitatively estimated.

\item A technical setup is given to arrange  the zeros of $B_{i}^{(r)}$ and $B_{i-2^{r-1}}^{(r-1)}B_{i+2^{r-1}}^{(r-1)}$ in pairs as shown in \eqref{theta} for the convenience of the forward error analysis. 
    Note that it is unclear how to arrange them in original paper \cite{Swarztrauber1974}.

\end{enumerate}

It is emphasized that  the main finding that the matrix polynomial $B_i^{(r)}$  corresponds to the eigenpolynomial of a principal submatrix of $-\mathcal{R}_n$ (see Theorem \ref{maintheorem} below) is the key ingredient of  the theoretical analysis of this work.

The rest of the paper is organized as followings. In section 2, cyclic reduction for solving the block-tridiagonal system is described.
In section 3, a brief description of the ECR algorithm is presented and the main  theorem is given. Section 4 presents  the forward error analysis  of the ECR algorithm in the reduction phase and the back-substitution phase.
\section{Cyclic reduction}

For  readers' convenience, the algorithm of cyclic reduction \cite{WalterGolub1997,BuzbeeGolub1970,Swarztrauber1977} is described here for solving block-tridiagonal systems that arise from discretizing the Poisson equation imposed on some rectangular domains by the finite difference method, which is of the form
\begin{align}\label{CRsystem}
\begin{bmatrix}
 A & I &     &  &  &  \\
 I & A & I &  &  &  \\
          & \ddots &\ddots &\ddots &\\
          &      &  I& A &I \\
         &      &         & I & A
\end{bmatrix}
\begin{bmatrix}
 u_1\\
 u_2\\
 \vdots\\
 u_n
\end{bmatrix}
=
\begin{bmatrix}
 g_1\\
 g_2\\
\vdots\\
 g_n
\end{bmatrix},
\end{align}
in which $u_i,g_i\in \RR^{m}$, $i=1,2,\cdots,n$, $I\in\RR^{m\times m}$ is the identity matrix with
 positive integers $m$ and $n$, and  the matrix $A\in\RR^{m\times m}$ is as follows
\begin{align*}
A=
\begin{bmatrix}
 -4 & 1 &     &  &    \\
 1 & -4 & 1 &  &    \\
          & \ddots &\ddots &\ddots &\\
    &      &  1& -4 & 1 \\
   &      &         & 1 & -4
\end{bmatrix}.
\end{align*}
The basic operation in cyclic reduction is the simultaneous elimination of  unknown vectors whose indices are odd. For the system \eqref{CRsystem}, the corresponding  elimination  can be done by matrix multiplications as follows.
Let $n=2^{k+1}-1$ and consider the following three consecutive equations
\begin{align*}
  u_{2j-2} + A u_{2j-1} + u_{2j}&=g_{2j-1},\\
  u_{2j-1} + A u_{2j} + u_{2j+1}& = g_{2j},\\
  u_{2j}+ A u_{2j+1} + u_{2j+2} &= g_{2j+1}.
\end{align*}
In order to eliminate $u_{2j-1}$ and $u_{2j+1}$, we multiply the second equation above with $-A$
and add all three equations. This leads to the following new equation
\begin{align}
  u_{2j-2} + (2I-A^2)u_{2j} + u_{2j+2}= g_{2j-1}-A \,g_{2j}+g_{2j+1}.
\end{align}
Let $A^{(0)}=A$ and $g_j^{(0)}=g_j$, and  define recursively
\begin{align}
  A^{(r+1)}&=2I-(A^{(r)})^2, \label{CR1}\\
  g_j^{(r+1)}&= g_{2j-1}^{(r)}-A^{(r)}g_{2j}^{(r)} + g_{2j+1}^{(r)}, \quad j=1,\cdots,2^{k+1-r}-1.\label{CR2}
\end{align}
Thus after $r$ reduction steps the remaining system of equations is of the size $(2^{k+1-r}-1)\times (2^{k+1-r}-1)$, which reads
\begin{align}\label{CR}
  u_{(j-1)2^r} + A^{(r)}u_{j2^r} + u_{(j+1)2^r}=g_j^{(r)}, \quad j=1,\cdots,2^{k+1-r}-1,
\end{align}
 here $u_0=u_{n+1}=0$. After $k$ reduction steps, the system with respect to one unknown vector $u_{2^k}$ is obtained as follows
\begin{align}\label{CRlastsystem}
  A^{(k)}u_{2^k}=g_1^{(k)},
\end{align}
where $A^{(k)}$ and $g_1^{(k)}$ are computed from \eqref{CR1} and \eqref{CR2}, respectively.
After determining $u_{2^k}$ a back-substitution is performed in which equation \eqref{CR} is
recursively solved for $u_{j2^r}$  while $u_{(j-1)2^r}$ and $u_{(j+1)2^r}$ are known from the  predecessor  level.

 In equation \eqref{CRlastsystem} and the back-substitution phase, the system of  equations with the matrices $A^{(r)}$ must be solved. Furthermore the transforming on the righthand side \eqref{CR2} needs matrix-vector multiplications with $A^{(r)}$. It follows from \eqref{CR1} that $A^{(r)}=P_{2^r}(A)$ is a polynomial of  matrix $A$ of degree $2^r$ and is  connected to the Chebyshev polynomials $T_n(t)$ as follows
\begin{align*}
   P_{2^r}(x)= -2T_{2^r}(-\frac{x}{2}), \quad r\geq 0.
\end{align*}
The zeros of $P_{2^r}(x)$ are as follows
\begin{align*}
 \lambda_i^{(r)}=-2\cos(\frac{2i-1}{2^{r+1}}\pi), \quad i=1,2,\cdots,2^r.
\end{align*}
Since for $r\geq 1$ the leading coefficient of $P_{2^r}(x)$ is $-1$, there holds the following factorization
\begin{align}\label{CRAr}
  P_{2^r}(x)=-\prod^{2^r}_{i=1}(x-\lambda_i^{(k)}).
\end{align}
Thus the matrices $A^{(r)}$ defined in equation \eqref{CR1} can be expressed as
\begin{align*}
 A^{(r)}= -\prod^{2^r}_{i=1}(A-\lambda_i^{(k)}I).
 \end{align*}

\section{The extended cyclic reduction algorithm }
The extended cyclic reduction (ECR) algorithm \cite{Swarztrauber1974} by Swarztrauber is used to solve the following more  general  block-tridiagonal system
\begin{align}\label{Glinearsystem}
\begin{bmatrix}
 B_1 & C_1 &     &  &    \\
 A_2 & B_2 & C_2 &  &    \\
          & \ddots &\ddots &\ddots & \\
          &      &  A_{n-1}& B_{n-1} &C_{n-1} \\
          &      &         & A_n & B_n
\end{bmatrix}
\begin{bmatrix}
 x_1\\
 x_2\\
 \vdots\\
 x_n
\end{bmatrix}
=
\begin{bmatrix}
 y_1\\
 y_2\\
\vdots\\
 y_n
\end{bmatrix},
\end{align}
where $n=2^k-1$ with some positive integer $k$, and $x_i,y_i\in \RR^{m}$, $i=1,2,\cdots,n$.
The matrices $A_i, B_i, C_i$ in \eqref{Glinearsystem} are of order $m$ and  of the form
\begin{align}\label{MatricesABC}
A_i = a_i I, \qquad B_i=B+b_i I,\qquad C_i = c_i I,
\end{align}
where $a_i,b_i$ and $c_i$ are scalars, and $B\in \RR^{m\times m}$ is a tridiagonal matrix. The system \eqref{Glinearsystem} can be rewritten in the tensor product formulation as follows
\begin{align}\label{system}
  \big( B\otimes I + I\otimes \mathcal{R}_n \big)X=Y,
\end{align}
where $X=(x_1^{\mathrm{ T }},x_2^{\mathrm{ T }},\ldots,x_n^{\mathrm{ T }})^{\mathrm{ T }}$, $Y=(y_1^{\mathrm{ T }},y_2^{\mathrm{ T }},\ldots,y_n^{\mathrm{ T }})^{\mathrm{ T }}$ and
\begin{align}\label{juzhenRn}
\mathcal{R}_n=
\begin{bmatrix}
 b_1 & c_1 &     &      \\
 a_2 & b_2 & c_2 &      \\
          & \ddots &\ddots &\ddots \\
     &      &  a_{n-1}& b_{n-1} &c_{n-1} \\
     &      &         & a_n & b_n
\end{bmatrix}.
\end{align}
\subsection{Brief overview and  implementation of the algorithm}
 For the sake of completeness,  a  brief overview and  some implementation issues
 of the ECR algorithm are described in this subsection.
We follow  the notation in \cite{Swarztrauber1974}.

The procedure begins with the reduction phase of the ECR algorithm.
In the reduction phase, the even  rows of the system \eqref{Glinearsystem} are eliminated at each step, and a reduced system with respect to the odd  unknowns is obtained. The size of the resulting system is about half  of the original system. After $\mathcal{O}(\log_2 n)$ steps we get one equation with one unknown vector, and after solving it there follows a back substitution phase during which the rest of the unknown vectors are computed. This is used as a basis for developing a fast algorithm.  The first system resulting from \eqref{Glinearsystem} is of block order $2^{k-1}-1$  with respect to the unknown vectors $x_2, x_4,\cdots, x_{2^k-2}$ by eliminating the unknown vectors $x_{i-1}$ and $x_{i+1}$ in the three block equations corresponding to block rows $i-1,i$ and $i+1$ of \eqref{Glinearsystem}.
Then one obtains the first system
\begin{align*}
   A_i^{(1)}x_{i-2}+  B_i^{(1)}x_i +  C_i^{(1)}x_{i+2}=y_i^{(1)},\quad i=2,4,\cdots,2^k-2,
\end{align*}
here
\begin{align*}
  A_i^{(1)}&=A_i B_{i+1}A_{i-1},\\
  B_i^{(1)}&=A_i B_{i+1}C_{i-1}-B_{i-1}B_{i+1}B_i+ C_i B_{i-1}A_{i+1},\\
  C_i^{(1)}&=C_i B_{i-1}C_{i+1},\\
  y_i^{(1)}&=A_i B_{i+1}y_{i-1}-B_{i-1}B_{i+1}y_i+ C_i B_{i-1}y_{i+1}.
\end{align*}
 This system is also  block tridiagonal and  the process above can be applied to  it. Then one obtains the next system  with respect to $2^{k-2}-1$ unknown vectors $x_4, x_8,\cdots, x_{2^k-4}$. What follows is the general reduction process which will now be described. Let $a_1=c_n=0$ and  for $i=1,2,\cdots,n$, define
\begin{align}\label{ABC}
  A_i^{(0)}=a_i I, \quad B_i^{(0)}= B+ b_i I, \quad C_i^{(0)}=c_i I, \quad y_i^{(0)}=y_i.
\end{align}
Let $\mathcal{D}_i^{(r+1)}=B_{i-2^{r-1}}^{(r-1)}B_{i+2^{r-1}}^{(r-1)}$,
for $r=0,1,\cdots,k-2$ and $i=2^{r+1}, 2\cdot 2^{r+1},\cdots, (2^{k-r-1})\cdot2^{r+1},$ define recursively
\begin{align}
  &A_i^{(r+1)}=(\mathcal{D}_i^{(r+1)})^{-1}\underbrace{ A_i^{(r)}B_{i+2^r}^{(r)}A_{i-2^r}^{(r)}}_{\widehat{A_i}^{(r+1)}},\label{Ahat}\\
  &B_i^{(r+1)}= (\mathcal{D}_i^{(r+1)})^{-1}\big(\underbrace{ A_i^{(r)}B_{i+2^r}^{(r)}C_{i-2^r}^{(r)}-B_{i-2^r}^{(r)}B_{i+2^r}^{(r)}B_{i}^{(r)}
  +C_i^{(r)}B_{i-2^r}^{(r)}A_{i+2^r}^{(r)}}_{\widehat{B_i}^{(r+1)}}\big),\label{Bhat}\\
  &C_i^{(r+1)}=(\mathcal{D}_i^{(r+1)})^{-1}\underbrace{ C_i^{(r)}B_{i-2^r}^{(r)}C_{i+2^r}^{(r)}}_{\widehat{C_i}^{(r+1)}},\label{Chat}\\
 & y_i^{(r+1)}= (\mathcal{D}_i^{(r+1)})^{-1}\big(A_i^{(r)}B_{i+2^r}^{(r)}y_{i-2^r}^{(r)}-B_{i-2^r}^{(r)}B_{i+2^r}^{(r)}y_{i}^{(r)}
  +B_{i-2^r}^{(r)}C_{i}^{(r)}y_{i+2^r}^{(r)}\big).
\end{align}
With $x_0=x_{2^{r+1}}=0,$ for $r=0,1,\cdots,k-2$, it yields a block tridiagonal system as follows
\begin{align}\label{kuaisanjiao}
  A_i^{(r)}x_{i-2^r}+ B_i^{(r)}x_i + C_i^{(r)}x_{i+2^r}= y_i^{(r)},\quad i=2^r, 2\cdot2^r,\cdots,(2^{k-r}-1)\cdot2^r.
\end{align}
Finally, for $r=k-1$, it arrives at the system only with respect to the unknown vector $x_{2^{k-1}}$,
\begin{align}\label{zuihouxitong}
  B^{(k-1)}_{2^{k-1}}x_{2^{k-1}}=y^{(k-1)}_{2^{k-1}}.
\end{align}
The back-substitution phase is initiated by solving \eqref{zuihouxitong} for unknown vector
$x_{2^{k-1}}$ and then proceeding backward using \eqref{kuaisanjiao}. For $r=k-2,k-3,\cdots,0$, and $i=2^r, 3\cdot2^r, 5\cdot2^r,\cdots, (2^{k-r}-1)\cdot 2^r,$ it leads to
\begin{align}\label{houdaixitong}
   x_i = (B_i^{(r)})^{-1}( y_i^{(r)} - A_i^{(r)}x_{i-2^r}- C_i^{(r)}x_{i+2^r} ).
\end{align}
The unknown vectors $x_{i-2^r}$ and $x_{i+2^r}$ on the right of \eqref{houdaixitong} are known from a predecessor step in the back-substitution.

It was proved in \cite{Swarztrauber1974} that $\mathcal{D}_i^{(r+1)}$ is a common factor of the matrices $A_i^{(r+1)}, B_i^{(r+1)}$ and $C_i^{(r+1)}$ defined in \eqref{Ahat}, \eqref{Bhat} and \eqref{Chat}, respectively. As a result, it follows from \eqref{ABC} and \eqref{Ahat} -- \eqref{Chat} that all the matrices $A_i^{(r)}, B_i^{(r)}, C_i^{(r)}$  are  polynomials of  matrix $B$.  Furthermore, the matrices $A_i^{(r)}$ and  $C_i^{(r)}$ can be expressed as follows (see \cite{Swarztrauber1974}):
\begin{align}
  A_i^{(r)}&= \alpha_i^{(r)}B_{i+2^{r-1}}^{(r-1)}, \quad \alpha_i^{(r)}=\prod_{j=i-2^{r}+1}^i a_j,\label{Ai}\\
  C_i^{(r)}&= \gamma_i^{(r)}B_{i-2^{r-1}}^{(r-1)}, \quad \gamma_i^{(r)}=\prod_{j=i}^{i+2^{r}-1}c_j,\label{Ci}
\end{align}
for $r=0,1,\cdots,k-1$ and $i=2^{r+1}, 2\cdot2^{r+1}, \cdots, (2^{k-r-1}-1)\cdot 2^{r+1}$, here
$B_i^{(-1)}=I$.
With the expression of $A_i^{(r)}$ in
\eqref{Ai} and $C_i^{(r)}$ in \eqref{Ci}, it is convenient to rewrite the ECR algorithm in terms of $B_i^{(r)}$,  which avoids the computation and storage of $A_i^{(r)}$ and $C_i^{(r)}$.
Hence, the preprocessing stage only consists  of computing the zeros of
 matrix polynomial $B_i^{(r)}$ with respect to $B$ as shown in \eqref{eq:phase0}.
This results in the ECR algorithm \cite{Swarztrauber1974} which contains the following   three phases:

\begin{enumerate}
\item[(0)]Preprocessing phase.  Instead of storing the matrices, it is to compute by  Newton's method and store the zeros of matrix polynomial $B_i^{(r)}$ that represent them, where  the matrix $B_i^{(r)}$ can be rewritten as follows
    \begin{align}\label{Bi0}
    B_i^{(-1)}=I,\quad  B_i^{(0)}=B_i, \qquad i=1,2,\ldots,n.
    \end{align}
and
\begin{align}\label{eq:phase0}
\begin{split}
 B_i^{(r)}=&\big( \mathcal{D}_i^{(r)}\big)^{-1}\big( \alpha_i^{(r-1)}\gamma_{i-2^{r-1}}^{(r-1)}B_{i+2^{r-2}}^{(r-2)}B_{i-3\cdot2^{r-2}}^{(r-2)}B_{i+2^{r-1}}^{(r-1)}\\
         &-B_{i-2^{r-1}}^{(r-1)}B_{i}^{(r-1)}B_{i+2^{r-1}}^{(r-1)}+ \alpha_{i+2^{r-1}}^{(r-1)}\gamma_{i}^{(r-1)}B_{i-2^{r-2}}^{(r-2)}B_{i+3\cdot2^{r-2}}^{(r-2)}B_{i-2^{r-1}}^{(r-1)} \big),
\end{split}
\end{align}
where  $\mathcal{D}_i^{(r)}=B_{i-2^{r-2}}^{(r-2)}B_{i+2^{r-2}}^{(r-2)}$ and
\begin{align*}
  & \alpha_1^{(0)}=0,\quad \alpha_i^{(0)}=a_i, \,i=2,3,\ldots,n,\quad \alpha_i^{(r)}= \alpha_i^{(r-1)}\alpha_{i-2^{r-1}}^{(r-1)},\\
  & \gamma_n^{(0)}=0,\quad \gamma_{i-1}^{(0)}=c_{i-1},\,
  i=2,3,\ldots,n,\quad \gamma_i^{(r)}= \gamma_i^{(r-1)}\gamma_{i+2^{r-1}}^{(r-1)},
\end{align*}
 for $r=1,2,\ldots,k-1$ and $i=2^{r},2\cdot2^{r},\ldots, (2^{k-r}-1)\cdot2^{r}$.

\item[(1)] Reduction phase.  Let $p_i^{(0)}=y_i$, $i=1,2,\ldots,n$. It is to compute $p_i^{(r+1)}$ for $r=0,1,\ldots,k-2$ and $i=2^{r+1},2\cdot2^{r+1},\ldots, (2^{k-r-1}-1)\cdot2^{r+1}$ by
\begin{align}\label{eq:phase1}
\begin{split}
  p_i^{(r+1)}= \alpha_i^{(r)}\big( B_{i-2^{r-1}}^{(r-1)} \big)^{-1}q_{i-2^r}^{(r)}+\gamma_i^{(r)}\big( B_{i+2^{r-1}}^{(r-1)}\big)^{-1}q_{i+2^r}^{(r)}-p_i^{(r)},
\end{split}
\end{align}
where
\begin{align}\label{qi}
  q_i^{(r)} = \big( B_i^{(r)}\big)^{-1}B_{i-2^{r-1}}^{(r-1)}B_{i+2^{r-1}}^{(r-1)}p_i^{(r)}.
\end{align}

\item[(2)]Back-substitution phase. With $x_0=x_{n+1}=0$, it is to compute $x_i$ for $r=k-1,k-2,\ldots,0$ and $i=2^{r},3\cdot2^r,\ldots,(2^{k-r}-1)\cdot2^{r}$ by
    \begin{align}\label{eq:phase2}
    \begin{split}
      x_i = \big( B_i^{(r)}\big)^{-1} B_{i-2^{r-1}}^{(r-1)}B_{i+2^{r-1}}^{(r-1)}\big[ p_i^{(r)}&-\alpha_i^{(r)}\big(B_{i-2^{r-1}}^{(r-1)}\big)^{-1}x_{i-2^r}\\
      &-\gamma_i^{(r)}\big(B_{i+2^{r-1}}^{(r-1)}\big)^{-1} x_{i+2^r}\big].
      \end{split}
    \end{align}
\end{enumerate}
There are two aspects for the implementation details of the ECR algorithm as follows:
\begin{enumerate}[1.]
\item  In the preprocessing phase, the matrix polynomial $B_i^{(r)}$ of degree $2^{r+1}-1$ is expressed in a factored form by using polynomials of degree one due to the fact that its  zeros computed by Newton's method are stored. In this way, the resulting method is satisfactory for high degree polynomials. However, it is difficult to get a quantitative  estimation of the computational cost and accuracy.  And the error is accumulated as the parameter $r$ increases, since the computation of the zeros of $B_i^{(r)}$ by Newton's method directly using \eqref{eq:phase0} depends on the zeros of $B_{i+2^{r-2}}^{(r-2)}, B_{i-3\cdot2^{r-2}}^{(r-2)}, B_{i+2^{r-1}}^{(r-1)}, B_{i-2^{r-1}}^{(r-1)}, B_{i}^{(r-1)}, B_{i-2^{r-1}}^{(r-2)}$ and $B_{i+3\cdot2^{r-2}}^{(r-2)}.$ 

\item In the reduction phase, to ensure a stable computing of vector
        $q_i^{(r)}$, the formulation in \eqref{qi} will not be directly used.
        It is necessary to multiply $p_i^{(r)}$ alternately by the inverse of a factor of  $B_i^{(r)}$ and by a factor of $B_{i-2^{r-1}}^{(r-1)}B_{i+2^{r-1}}^{(r-1)}$. That is, $q_i^{(r)}$ is obtained by first defining $z_0=p_i^{(r)}$ and computing $z_j$ recursively by solving the following linear systems
       \begin{align}\label{thetaphi}
         (B-\theta_j I)z_{j+1}=(B-\phi_j I)z_{j},\quad j=1,2,\cdots,2^{r+1}-2,
       \end{align}
       where $\theta_j$ is a zero of $B_i^{(r)}$ and $\phi_j$ is a zero of $B_{i-2^{r-1}}^{(r-1)}B_{i+2^{r-1}}^{(r-1)}$. Then $q_i^{(r)}$ is given by
      $(B-\theta_{2^{r+1}-1} I)^{-1}z_{2^{r+1}-1}$. The implementation of \eqref{thetaphi} is carried out through defining $\delta_{j+1}$ by
     \begin{align}\label{zj+1}
       z_{j+1}=\delta_{j+1}+z_j.
     \end{align}
     Then substituting \eqref{zj+1} into \eqref{thetaphi}, one obtains
     \begin{align}\label{zz}
       ( B-\theta_j I)\delta_{j+1}=(\theta_j-\phi_j)z_j.
     \end{align}
      After $\delta_{j+1}$ has been obtained from \eqref{zz}, $z_{j+1}$ can be got from \eqref{zj+1}. In this way, all matrix multiplications in \eqref{qi} can be avoided. Here $\theta_j$ and $\phi_j$ are selected so that $\max_{j}|\theta_j-\phi_j|$ is as small as possible so that  the roundoff error will not grow.
       Such a  technique  is also  used to compute a term like $\alpha_i^{(r)}\big( B_{i-2^{r-1}}^{(r-1)} \big)^{-1}q_{i-2^r}^{(r)}$ in \eqref{eq:phase1}, since it is not possible to avoid repeated multiplications by the inverse of factors of $B_{i-2^{r-1}}^{(r-1)}$  which results in error.   The technique   utilized in the reduction phase should also be used in the implementation of the back-substitution phase in \eqref{eq:phase2}.
It is   noteworthy that it is unclear how to select $\theta_j$ and $\phi_j$  from the paper \cite{Swarztrauber1974}.

\end{enumerate}

The improvements we make corresponding to the problems above are as follows:
\begin{enumerate}[1.]
 \item Thanks to the finding that the zeros of  matrix polynomial $B_i^{(r)}$ are the eigenvalues of a principal submatrix of $-\mathcal{R}_n$, the MOB method for the eigenvalues of a symmetric tridiagonal matrix is applied  to compute the  zeros of  $B_{i}^{(r)}$,  the computational cost and accuracy of which are  quantitatively estimated. And the accumulation of the error will not occur since the recurrence relation \eqref{eq:phase0} is not used.

 \item In the forthcoming theoretical analysis (see Theorem \ref{dingli1} and Theorem \ref{dingli2}), we choose $\theta_j$ and $\phi_j$ in pairs in order after arranging the zeros from the largest one to the smallest one, that is,
     \begin{align}\label{theta}
   \theta_j=\mu_{j+1}, \quad \phi_j= \lambda_j, \quad j=1,2,\cdots, 2^{r+1}-2,
    \end{align}
where $\{ \mu_i \}_{i=1}^{2^{r+1}-1}$ are the zeros of matrix polynomial  $B_{i}^{(r)} $ such that $0> \mu_1> \mu_2 >\ldots >\mu_{2^{r+1}-1}$, and $\{ \lambda_i \}_{i=1}^{2^{r+1}-2}$ are the zeros of the corresponding polynomial of matrix  $B_{i-2^{r-1}}^{(r-1)}B_{i+2^{r-1}}^{(r-1)}$  such that $0> \lambda_1\geq \lambda_2 \geq \ldots \geq \lambda_{2^{r+1}-2}$.  It is proved in
Theorem \ref{dingli1} and Theorem \ref{dingli2}  that the algorithm with such an arrangement is stable.  
\end{enumerate}

\subsection{ Results on the matrix polynomial $B_i^{(r)}$  }
This subsection will prove the  matrix polynomial
  $B_i^{(r)}$ corresponds to the eigenpolynomial of a principal submatrix of $-\mathcal{R}_n$.
To this end, we first introduce the following notation.

For $\mathcal{R}_n\in \RR^{n\times n}$ in \eqref{juzhenRn} and $2\le k\le n$, define
\begin{align*}
&\mathcal{F}_{k}:= \mathcal{R}_n([1:k;1:k]),\\
 &\mathcal{L}_{k}:= \mathcal{R}_n([2:k;2:k]),
\end{align*}
 that is, $\mathcal{F}_{k}$ is  the matrix formed from the first $k$ rows and $k$ columns of $\mathcal{R}_n$, and
$\mathcal{L}_{k}$ is  the matrix formed from the second row to the $k$th row and the second column to the $k$th column of $\mathcal{R}_n$.
\begin{lemma}\label{lemma}
It holds that
\begin{align}\label{lem:jielun}
  \det(\mathcal{F}_{n-1})\det(\mathcal{L}_{n})-\det(\mathcal{F}_n)\det(\mathcal{L}_{n-1})= \prod_{i=1}^{n-1}a_{i+1}c_i.
\end{align}
\end{lemma}

\begin{proof}
 The proof will now proceed by induction on the order of matrix $\mathcal{R}_n$. If $n=3$, the direct computation leads to
 \begin{align*}
   \det\Big(\begin{bmatrix}
 b_1 & c_1   \\
 a_2 & b_2
\end{bmatrix}\Big)\cdot
\det\Big(\begin{bmatrix}
 b_2 & c_2   \\
 a_3 & b_3
\end{bmatrix}\Big)-b_2\cdot
\det\Big(\begin{bmatrix}
 b_1 & c_1 & 0  \\
 a_2 & b_2 & c_2\\
 0 &   a_3 & b_3
\end{bmatrix}\Big)
=a_2 c_1 a_3 c_2.
 \end{align*}

 Assume \eqref{lem:jielun} holds for $n=t-1$, that is
 \begin{align*}
  \det(\mathcal{F}_{t-2})\det(\mathcal{L}_{t-1})-\det(\mathcal{F}_{t-1})\det(\mathcal{L}_{t-2})= \prod\limits_{i=1}^{t-2}a_{i+1}c_i.
\end{align*}
Let
\begin{align*}
 & \widehat{\mathbf{c}}_1=(c_1,0,\ldots,0)\in\RR^{t-2},\,\,\quad
 \widehat{\mathbf{c}}_{t-1}=(0,\ldots,0,c_{t-1})^{\mathrm{ T }}\in\RR^{t-2},\\
 &\widehat{\mathbf{a}}_t=(0,\ldots,0,a_t)\in\RR^{t-2}, \quad
 \widehat{\mathbf{a}}_2=(a_2,0,\ldots,0)^{\mathrm{ T }}\in\RR^{t-2},
\end{align*}
and
\begin{align*}
  &\mathcal{A}_{1}=
\begin{bmatrix}
 b_1                   & \mathbf{\widehat{c}_1}  \\
 \mathbf{\widehat{a}_2} & \mathcal{L}_{t-1}\\
\end{bmatrix},
\mathcal{A}_{4}=
\begin{bmatrix}
 \mathcal{L}_{t-1} & \mathbf{\widehat{c}_{t-1}}\\
  \mathbf{\widehat{a}_t} & b_t
\end{bmatrix},
  \mathcal{A}_{5}=
  \begin{bmatrix}
     b_t        & 0\\
    0& \mathcal{L}_{t-1}
   \end{bmatrix},\quad\\
 &\mathcal{A}_{3}=
 \begin{bmatrix}
  0&  0\\
 \mathbf{\widehat{c}_{t-1}}& 0
 \end{bmatrix},\quad
 \mathcal{A}_{22}=
\begin{bmatrix}
 0 & 0 \\
 0 & \mathbf{\widehat{a}_t}
\end{bmatrix},\quad
\mathcal{A}_{33}=
 \begin{bmatrix}
  0&  0\\
   0& \mathbf{\widehat{c}_{t-1}}
 \end{bmatrix},\quad\\
& \mathcal{A}_{45}=
\begin{bmatrix}
 \mathcal{L}_{t-1} & 0\\
  \mathbf{\widehat{a}}_t & b_t
\end{bmatrix},\quad
\mathcal{A}_{55}=
\begin{bmatrix}
     \mathcal{L}_{t-1} & 0\\
    0& b_t
   \end{bmatrix},\quad
\mathcal{A}_{2}=
\begin{bmatrix}
 0 &  \mathbf{\widehat{a}_t}\\
 0 & 0
\end{bmatrix}.
\end{align*}
For $n=t$, it follows that
\begin{align*}
&\det(\mathcal{F}_{t-1})\det(\mathcal{L}_{t})-\det(\mathcal{F}_{t})\det(\mathcal{L}_{t-1})
=\det\Big(
  \begin{bmatrix}
        \mathcal{A}_{1} & 0 \\
           0 & \mathcal{A}_{4} \\
\end{bmatrix}\Big)
-
\det\Big(
  \begin{bmatrix}
        \mathcal{A}_{1} & \mathcal{A}_{3} \\
         \mathcal{A}_{2} & \mathcal{A}_{5} \\
\end{bmatrix}\Big)\\
= &\det\Big(
  \begin{bmatrix}
        \mathcal{A}_{1} & 0 \\
           0 & \mathcal{A}_{4} \\
\end{bmatrix}\Big)
-
\det\Big(
  \begin{bmatrix}
        \mathcal{A}_{1} & \mathcal{A}_{33} \\
           0& \mathcal{A}_{45} \\
\end{bmatrix}\Big)
+
\det\Big(
  \begin{bmatrix}
        \mathcal{A}_{1} & \mathcal{A}_{33} \\
           0 & \mathcal{A}_{45} \\
\end{bmatrix}\Big)
-
\det\Big(
  \begin{bmatrix}
        \mathcal{A}_{1} & \mathcal{A}_{33} \\
         \mathcal{A}_{22} & \mathcal{A}_{55} \\
\end{bmatrix}\Big).
\end{align*}
Moreover, it gets from the property of the determinant of matrices that
\begin{align*}
& \det\Big(
  \begin{bmatrix}
        \mathcal{A}_{1} & 0 \\
           0 & \mathcal{A}_{4} \\
\end{bmatrix}\Big)
-
\det\Big(
  \begin{bmatrix}
        \mathcal{A}_{1} & \mathcal{A}_{33} \\
           0& \mathcal{A}_{45} \\
\end{bmatrix}\Big)\\
=& \det\Big(
\begin{bmatrix}
 b_1                   & \mathbf{\widehat{c}_1}      & 0 & 0 \\
 \mathbf{\widehat{a}_2} & \mathcal{L}_{t-1} & 0      & -\mathbf{\widehat{c}_{t-1}}\\
 0 &  0 & \mathcal{L}_{t-1} & \mathbf{\widehat{c}_{t-1}}\\
 0 & 0  & \mathbf{\widehat{a}_t} & 0
\end{bmatrix}\Big)
=-a_t c_{t-1}\det\Big( \mathcal{L}_{t-2}\Big)\det\Big( \mathcal{F}_{t-1}\Big),
\end{align*}
Similarly, it holds that
\begin{align*}
 & \det\Big(
  \begin{bmatrix}
        \mathcal{A}_{1} & \mathcal{A}_{33} \\
           0 & \mathcal{A}_{45} \\
\end{bmatrix}\Big)
-
\det\Big(
  \begin{bmatrix}
        \mathcal{A}_{1} & \mathcal{A}_{33} \\
         \mathcal{A}_{22} & \mathcal{A}_{55} \\
\end{bmatrix}\Big)\\
= &\det\Big(
  \begin{bmatrix}
 b_1                   & \mathbf{\widehat{c}_1}      & 0 & 0 \\
 \mathbf{\widehat{a}_2} & \mathcal{L}_{t-1} & 0     & \mathbf{\widehat{c}_{t-1}} \\
 0 &  0       & \mathcal{L}_{t-1}        & 0\\
 0 &  -\mathbf{\widehat{a}_t} & \mathbf{\widehat{a}_t}& 0
\end{bmatrix}\Big)
=a_tc_{t-1}\det\Big( \mathcal{L}_{t-1} \Big)\det\Big( \mathcal{F}_{t-2} \Big).
\end{align*}
 Since  $\mathcal{F}_{t-1}\in \RR^{(t-1)\times (t-1)}$, it follows from the induction assumption that
\begin{align*}
  \det\Big( \mathcal{L}_{t-1} \Big)\det\Big( \mathcal{F}_{t-2}\Big)-\det\Big( \mathcal{L}_{t-2}\Big)\det\Big( \mathcal{F}_{t-1}\Big)
  = \prod\limits_{i=1}^{t-2}a_{i+1}c_i.
\end{align*}
This leads to
\begin{align*}
   \det(\mathcal{F}_{t-1})\det(\mathcal{L}_{t})-\det(\mathcal{F}_{t})\det(\mathcal{L}_{t-1})
   =a_t c_{t-1}\prod\limits_{i=1}^{t-2}a_{i+1}c_i
   =\prod\limits_{i=1}^{t-1}a_{i+1}c_i,
\end{align*}
which proves the conclusion \eqref{lem:jielun}.
\end{proof}

Lemma \ref{lemma} plays an important role in the proof of the theorem below.


\begin{theorem}\label{maintheorem}
Define
\begin{align*}
\mathcal{R}_i^{(r)}:= -\mathcal{R}_n\Big(\big[i-(2^r-1):i+(2^r-1);i-(2^r-1):i+(2^r-1)\big]\Big),
\end{align*}
for $r=0,1,\cdots,k-1$ and $i=2^r, 2\cdot2^r,\cdots,(2^{k-r}-1)\cdot2^r.$
Let $f_i^{(r)}(x)$ be the eigenpolynomial of matrix $\mathcal{R}_i^{(r)}$, i.e., $$f_i^{(r)}(x)=\det\left(xI-\mathcal{R}_i^{(r)}\right).$$ Then it holds that
\begin{align}\label{mainresult}
 B_i^{(r)} = (-1)^rf_i^{(r)}(B).
\end{align}
\end{theorem}

\begin{proof}
 We prove the theorem  by induction on $r$.  In view of \eqref{Bi0} and
 $\mathcal{R}_i^{(0)} = -\big[ b_i \big],\, i=1,2,\cdots,n,$  it is obvious that \eqref{mainresult} holds for $r = 0$.
For $r=1$,
\begin{align*}
  \mathcal{R}_i^{(1)}=-\begin{bmatrix} b_{i-1} & c_{i-1}&0\\
a_{i}& b_{i} & c_{i}\\ 0 & a_{i+1} & b_{i+1}  \end{bmatrix}, \quad i=2,4,\cdots, n-1.
\end{align*}
It follows  that
\begin{align*}
&f_i^{(1)}(x)=\det\left(xI-\mathcal{R}_i^{(1)}\right)= \det\Big( \begin{bmatrix} x+b_{i-1} & c_{i-1}&0\\
a_{i}& x+b_{i} & c_{i}\\ 0 & a_{i+1} & x+b_{i+1}  \end{bmatrix}\Big)\\
= &(x+b_{i-1})(x+b_{i+1})(x+b_{i}) - a_{i+1}c_i(x+b_{i-1}) - a_ic_{i-1}(x+b_{i+1}),
\end{align*}
which indicates  $B_i^{(1)} = -f_i^{(1)}(B)$,  where $B_i^{(1)}$ is the  polynomial of tridiagonal matrix $B$ as shown in \eqref{eq:phase0} with $r=1$.
Assume \eqref{mainresult} holds for $r\leq t-1$. Next, we turn to the case  $r=t.$

In view of \eqref{eq:phase0}, it holds the following identity
\begin{align}\label{eq:det}
\begin{split}
\mathcal{L}:=&\det\left(xI- \mathcal{R}_i^{(t)}\right)\det\left(xI - \mathcal{R}_{i-2^{t-2}}^{(t-2)}\right)\det\left(xI - \mathcal{R}_{i+2^{t-2}}^{(t-2)}\right)\\
 =& \det\Big(xI -  \begin{bmatrix}\mathcal{R}_{i-2^{t-1}}^{(t-1)} & 0 &-\mathbf{\widehat{c}}_{i-1} & & \\
0 & \mathcal{R}_{i-2^{t-2}}^{(t-2)} & 0 & & \\
-\mathbf{\widehat{a}}_i & 0 & -b_i & 0 & -\mathbf{\widehat{c}}_{i} \\
 & & 0 & \mathcal{R}_{i+2^{t-2}}^{(t-2)}  & 0 \\
 & & -\mathbf{\widehat{a}}_{i+1} & 0& \mathcal{R}_{i+2^{t-1}}^{(t-1)} \\
 \end{bmatrix}\Big)\\
=& f_{i-2^{t-1}}^{(t-1)}(x)\,f_{i+2^{t-1}}^{(t-1)}(x)
  \det\Big(xI -  \begin{bmatrix}
 \mathcal{R}_{i-2^{t-2}}^{(t-2)} & 0 & \\
 0 & -b_i & 0  \\
    & 0 & \mathcal{R}_{i+2^{t-2}}^{(t-2)}  \\
 \end{bmatrix}\Big)\\
&-a_{i+1}c_i f_{i-2^{t-1}}^{(t-1)}(x)\,f_{i-2^{t-2}}^{(t-2)}(x)\, f_{i+2^{t-2}}^{(t-2)}(x)\,
\det \Big( xI- \widetilde{\mathcal{R}}_{i+2^{t-1}}^{(t-1)}\Big)\\
& -a_ic_{i-1} f_{i+2^{t-1}}^{(t-1)}(x)\,f_{i-2^{t-2}}^{(t-2)}(x)\,f_{i+2^{t-2}}^{(t-2)}(x)\,
\det\Big(xI - \widehat{\mathcal{R}}_{i-2^{t-1}}^{(t-1)}\Big),
\end{split}
\end{align}
 where
\begin{align*}
&\mathbf{\widehat{a}}_{i} = (0,\cdots,0,a_i)\in \RR^{2^t-1}, \quad \mathbf{\widehat{c}}_{i-1} =  (0,\cdots,0,c_{i-1})^{\mathrm{ T }}\in\RR^{2^t-1},\\
&\mathbf{\widehat{c}}_{i} =  (c_i,0,\cdots,0)\in \RR^{2^t-1},\quad
\mathbf{\widehat{a}}_{i+1} = (a_{i+1},0,\cdots,0)^{\mathrm{ T }}\in \RR^{2^t-1},\\
& \widehat{\mathcal{R}}_{i-2^{t-1}}^{(t-1)}:= -\mathcal{R}_{n}\big([i-2^t+1:i-2;i-2^t+1:i-2]\big),\\
& \widetilde{\mathcal{R}}_{i+2^{t-1}}^{(t-1)}:= -\mathcal{R}_{n}\big([i+2:i+2^t-1;i+2:i+2^t-1]\big).
\end{align*}
Note that
 \begin{align}\label{jielun}
 \begin{split}
 \det\Big(xI - \begin{bmatrix}
 \mathcal{R}_{i-2^{t-2}}^{(t-2)} & 0 & \\
0 & -b_i & 0  \\
 & 0 & \mathcal{R}_{i+2^{t-2}}^{(t-2)}  \\
 \end{bmatrix}&\Big)
 = f_{i}^{(t-1)}(x)
 + a_{i}c_{i-1}\det\left(xI- \widehat{\mathcal{R}}_{i-2^{t-2}}^{(t-2)}\right)f_{i+2^{t-2}}^{(t-2)}(x)\\
  &+a_{i+1}c_{i}\det\left(xI- \widetilde{\mathcal{R}}_{i+2^{t-2}}^{(t-2)}  \right) f_{i-2^{t-2}}^{(t-2)}(x).
\end{split}
\end{align}
A substitution of \eqref{jielun} into \eqref{eq:det} leads to
\begin{align*}
  \mathcal{L}&= f_{i-2^{t-1}}^{(t-1)}(x)f_{i+2^{t-1}}^{(t-1)}(x)f_{i}^{(t-1)}(x)\\
  &+a_ic_{i-1}\, f_{i+2^{t-1}}^{(t-1)}(x)f_{i+2^{t-2}}^{(t-2)}(x)
  \Big(f_{i-2^{t-1}}^{(t-1)}\det\Big(xI-\widehat{\mathcal{R}}_{i-2^{t-2}}^{(t-2)} \Big)
 -f_{i-2^{t-2}}^{(t-2)}\det\Big(xI - \widehat{\mathcal{R}}_{i-2^{t-1}}^{(t-1)}\Big)\Big)\\
  &+a_{i+1}c_i\, f_{i-2^{t-1}}^{(t-1)}(x)f_{i-2^{t-2}}^{(t-2)}(x)
  \Big( f_{i+2^{t-1}}^{(t-1)}\det\Big(xI - \widetilde{\mathcal{R}}_{i+2^{t-2}}^{(t-2)}\Big)-f_{i+2^{t-2}}^{(t-2)}
\det \Big( xI- \widetilde{\mathcal{R}}_{i+2^{t-1}}^{(t-1)}\Big)\Big).
\end{align*}
Furthermore, it will be verified that
\begin{align}\label{gongshi1}
\begin{split}
&f_{i-2^{t-1}}^{(t-1)}
\det\Big( xI\! - \! \widehat{\mathcal{R}}_{i-2^{t-2}}^{(t-2)}\Big)
\!-\! f_{i-2^{t-2}}^{(t-2)}\det\Big(xI - \widehat{\mathcal{R}}_{i-2^{t-1}}^{(t-1)}\Big)
\!=\!  -f_{i-3\cdot 2^{t-2}}^{(t-2)}\prod\limits_{j=i-2^{t-1}}^{i-2} a_{j+1}c_j,
\end{split}
\end{align}
\begin{align}\label{gongshi2}
\begin{split}
&f_{i+2^{t-1}}^{(t-1)}
\det\left( xI \!-\!  \widetilde{\mathcal{R}}_{i+2^{t-2}}^{(t-2)}\right)
\!-\! f_{i+2^{t-2}}^{(t-2)}\det\left(xI - \widetilde{\mathcal{R}}_{i+2^{t-1}}^{(t-1)}\right)
\!=\! -f_{i+3\cdot 2^{t-2}}^{(t-2)}
\prod\limits^{i+2^{t-1}-1}_{j=i+1} a_{j+1}c_j.
\end{split}
\end{align}
The proof of \eqref{gongshi1} and \eqref{gongshi2} is postponed to Appendix A.
Finally, it follows from  \eqref{gongshi1} and \eqref{gongshi2} that
\begin{align*}
\mathcal{L}&=f_{i-2^{t-1}}^{(t-1)}(x)f_{i+2^{t-1}}^{(t-1)}(x)f_{i}^{(t-1)}(x)
   - f_{i+2^{t-1}}^{(t-1)}(x)f_{i+2^{t-2}}^{(t-2)}(x)
 f_{i-3\cdot 2^{t-2}}^{(t-2)}(x)\prod\limits_{j=i-2^{t-1}}^{i-1} a_{j+1}c_j\\
  &-f_{i-2^{t-1}}^{(t-1)}(x) f_{i-2^{t-2}}^{(t-2)}(x)f_{i+3\cdot 2^{t-2}}^{(t-2)}(x)
\prod\limits^{i+2^{t-1}-1}_{j=i} a_{j+1}c_j,
\end{align*}
which is the corresponding polynomial of tridiagonal  matrix $B$ in \eqref{eq:phase0} with $r=t$. The proof is completed.
\end{proof}

\begin{remark}
  Theorem \ref{maintheorem}  indicates that $B_i^{(r)}$ is  naturally a  polynomial of the tridiagonal matrix $B$ and $\widehat{B_i}^{(r)}$ shown in \eqref{Bhat} contains the factor $B_{i-2^{r-2}}^{(r-2)}B_{i+2^{r-2}}^{(r-2)}$.
\end{remark}

\section{ Forward error  analysis of the extended cyclic reduction algorithm }

 In this section, forward error analysis of the  ECR algorithm for the linear algebraic system \eqref{Glinearsystem} is presented.

 \subsection{ Zeros of matrix polynomial $B^{(r)}_{i}$ }
 In \cite{Swarztrauber1974}, the zeros of  matrix polynomial $B_i^{(r)}$ of matrix $B$ are computed
 by  Newton's method. It is observed that Newton's method therein  is globally convergent only if the zeros are real and simple from numerical experiments. However, a theoretical proof is missing in literatures. In what follows, a  theoretical  analysis is given to show that  the zeros are real and simple.

 A combination of Theorem \ref{maintheorem} and the following lemma arrives at the conclusion about the zeros of  matrix polynomial $B_i^{(r)}$ .
\begin{lemma}[\cite{Stewart1990}]\label{lem:eigen}
  Let $H\in\mathbb R^{n\times n}$ be the arrowhead matrix of the form
\begin{equation}
H = \begin{bmatrix}
\alpha & \mathbf{z}^{\mathrm{ T }} \\
\mathbf{z} & D
\end{bmatrix},
\end{equation}
where
\begin{align*}
&\mathbf{z}=(z_2,z_3,\cdots,z_n)^{\mathrm{ T }},\qquad z_j\neq 0,\,\, j=2,3,\cdots,n, \\
& D = \mathrm{diag}(d_2,d_3,\cdots,d_n), \quad d_2<d_3<\cdots < d_n.
\end{align*}
Assume the eigenvalues of $H$ are in the order $\lambda_1\leq \lambda_2\leq \cdots\leq \lambda_n$. Then it holds that
$$
\lambda_1<d_2<\lambda_2<d_3<\cdots<\lambda_{n-1} < d_n <\lambda_n.
$$
\end{lemma}

\begin{theorem}
 The zeros of matrix polynomial $B_i^{(r)}$ are real and simple.
\end{theorem}
\begin{proof}
  Due to the fact that $\mathcal{R}_i^{(r)}$ is a real symmetric negative definite matrix, it yields that the zeros of $B_i^{(r)}$ are real.
 In what follows, we prove the eigenvalues of any principal submatrix of $\mathcal{R}_n$ are simple based on an induction argument.  This shows  that the zeros of $B_i^{(r)}$ are simple.

 It is obvious that the eigenvalue of any principal submatrix of order one of $\mathcal{R}_n$ is simple. Any principal submatrix  of order two of form
$
   \begin{bmatrix}
      b_i & a_{i+1} \\
      a_{i+1} & b_{i+1}
\end{bmatrix}
$
  has two eigenvalues $\lambda_1$ and $\lambda_2$($\lambda_1\leq \lambda_2$) satisfying $\lambda_1<b_{i+1}<\lambda_2$ due to  Lemma \ref{lem:eigen}. Assume eigenvalues of
   any principal submatrix of order $t-1$ of $\mathcal{R}_n$ are simple, then eigenvalues of principal submatrix of order $t$ will be investigated.

   Denote a principal submatrix of order $t$ of $\mathcal{R}_n$ by
  \begin{align*}
  T_t=
\begin{bmatrix}
b_{i} & a_{i+1} & & &  \\
a_{i+1} & b_{i+1} & a_{i+2}  \\
&\ddots & \ddots & \ddots \\
 & & a_{i+t-2} & b_{i+t-2} & a_{i+t-1} \\
  & & & a_{i+t-1} & b_{i+t-1}
\end{bmatrix}
 ,\quad i\geq 1,\quad i+t \leq n.
\end{align*}
Note that $T_t = \begin{bmatrix}
b_i & \widetilde{\mathbf{a}}_{i+1}^{\mathrm{T}} \\
\widetilde{\mathbf{a}}_{i+1} & T_{t-1}
\end{bmatrix},$
where $\widetilde{\mathbf{a}}_{i+1} = (a_{i+1},0,\cdots,0)^{\mathrm{ T }}\in\RR^{t-1}$, $T_{t-1}$ is the  principal submatrix of order $t-1$ of $\mathcal{R}_n$. Since $T_{t-1}$ is a symmetric matrix, there exists the following orthogonal decomposition
\begin{align}\label{eq:pufenjie}
  T_{t-1}V=VD,
\end{align}
where $V$ is an orthogonal matrix, $D=\text{diag}(d_2,d_3,\ldots,d_t)$. It follows from the inductive hypothesis that there are $t-1$ distinct eigenvalues for matrix $T_{t-1}$  and $0<d_2<d_3<\ldots <d_t$. Let $v_j$ be the $j$-th column of matrix $V$, it holds that
\begin{align*}
  T_{t-1}v_j = d_j v_j,  \qquad  j=2,3,\ldots,t.
\end{align*}
It is straightforward to see that the first component of vector $v_j$ is non-zero. Then  the first row of matrix $V$ is also non-zero. Thanks to \eqref{eq:pufenjie}, it holds that
 \begin{align*}
T_t = \begin{bmatrix} 1&0\\ 0&V \end{bmatrix}
\begin{bmatrix} b_i & \widetilde{\mathbf{a}}_{i+1}^{\mathrm{T}} V\\ V^{\mathrm{T}} \widetilde{\mathbf{a}}_{i+1} & D \end{bmatrix}
\begin{bmatrix} 1&0\\ 0&V^{\mathrm{T}} \end{bmatrix}.
\end{align*}
Note that all components of vector $\widetilde{\mathbf{a}}_{i+1}^{\mathrm{T}} V$ are non-zero, it gets from Lemma \ref{lem:eigen} that the eigenvalues of  matrix $T_t$ are simple. Finally, this implies that the zeros of  $B_i^{(r)}$ are simple.
\end{proof}

 Thanks to the fact that the zeros of matrix polynomial $B_{i}^{(r)}$ correspond to the eigenvalues of a principal submatrix of $-\mathcal{R}_n$, the method of bisection (MOB) \cite{Barth1967} for the eigenvalues of a symmetric tridiagonal matrix is applied  to compute the  zeros of  $B_{i}^{(r)}$, the computation complexity  of which is of order $\mathcal{O}(n^2)$ for computing all the eigenvalues of a matrix of order $n$. In addition, the following error estimate was proved in \cite{Barth1967},
  \begin{align*}
    \varepsilon=\max|\lambda_j-\widetilde{\lambda}_j|\le \frac{15}{2} \kappa\cdot
    \max_{1\le i\le n}|b_i\pm (|a_i|+|a_{i+1}|) |, \quad a_1=0, a_{n+1}=0,
  \end{align*}
  where $\lambda_j$ and $\widetilde{\lambda}_j$ are the $j$th exact and approximate eigenvalues
  by the MOB, respectively, $b_i$  and $a_i$ are  diagonal elements and
subdiagonal elements ($a_i\neq 0$ for $i\neq 1$) of the symmetric tridiagonal matrix of order $n$, respectively,
 $\kappa$ is a preassigned tolerance and is set to be the unit roundoff $\mathbf{u}$ in this paper.
  In  what follows,  the zeros of  $B_{i}^{(r)}$  and $B_{i-2^{r-1}}^{(r-1)}B_{i+2^{r-1}}^{(r-1)}$ are investigated.

\begin{lemma}\label{lem:4.2}
  Let $\{ \mu_i \}_{i=1}^{2^{r+1}-1}$ be the zeros of  $B_{i}^{(r)}$ such that $0> \mu_1> \mu_2 > \ldots > \mu_{2^{r+1}-1}$, and
  $\{ \lambda_i \}_{i=1}^{2^{r+1}-2}$ be the zeros of $B_{i-2^{r-1}}^{(r-1)} B_{i+2^{r-1}}^{(r-1)}$  such that $0> \lambda_1\geq \lambda_2 \geq \ldots \geq \lambda_{2^{r+1}-2}$, it holds that
 \begin{align*}
   \frac{ |\mu_{\ell+1}-\lambda_{\ell}| }{ |\mu_{\ell+1}| }+ \frac{|\mu_{\ell}|}{|\mu_{\ell+1}|}  <1, \qquad  \ell=1,2,\ldots,2^{r+1}-2.
 \end{align*}
\end{lemma}

\begin{proof}
 It follows from the Sturm sequence property  \cite{GolubVan2013} that
 $ |\lambda_{\ell}| > |\mu_{\ell}|$. Then it yields that
 \begin{align*}
   1= \frac{ |\mu_{\ell+1}-\lambda_{\ell}| }{ |\mu_{\ell+1}| }+ \frac{|\lambda_{\ell}|}{|\mu_{\ell+1}|}
   >  \frac{ |\mu_{\ell+1}-\lambda_{\ell}| }{ |\mu_{\ell+1}| }+ \frac{|\mu_{\ell}|}{|\mu_{\ell+1}|},\quad \ell=1,2,\ldots,2^{r+1}-2,
 \end{align*}
 which ends the proof.
\end{proof}

\subsection{ Forward error analysis of  Gaussian elimination for solving a tridiagonal system}
In this subsection, the forward error analysis of  Gaussian elimination for solving a tridiagonal system is given, since the ECR algorithm is based on the computation of linear systems with the symmetric tridiagonal matrix.

The following componentwise backward error analysis result is useful in the round-off error analysis. Throughout  this paper,  the unit roundoff $\mathbf{u}$ is assumed to be sufficiently small.
\begin{lemma}[\cite{Higham1990},Theorem 3.2]\label{lem1}
 If a tridiagonal matrix $A$ is  symmetric positive definite, then Gaussian elimination for solving $Ax=b$ succeeds and the computed solution $\widehat{x}$ satisfies
 \begin{align}
  &(A+G)\widehat{x}=b, \label{eq1}\\
  & |G|\leq g(\mathbf{u})|A|, \quad g(\mathbf{u}):=\frac{4\mathbf{u}+3\mathbf{u}^2+ \mathbf{u}^3}{1-\mathbf{u}}\label{eq2},
 \end{align}
where  the backward error matrix $G$ is small componentwise relative to $A$,
 the absolute value operation $|\cdot|$ and  the matrix inequality are interpreted componentwise.
\end{lemma}

Applying the standard perturbation theory to \eqref{eq1}, one obtains the following forward error bound.
\begin{lemma}\label{lem:wucha}
  Let $A$ be a   symmetric positive definite tridiagonal matrix. Assume that the computed solution $\widehat{x}$  of  the system $Ax=b$  by Gaussian elimination satisfies
  \begin{align}\label{fangcheng}
   (A+G)\widehat{x}=b.
  \end{align}
  Then, it holds that
    \begin{align}\label{jielun1}
    \|x-\widehat{x}\|_2 \le \xi(A,\mathbf{u})\|x\|_2, \quad  \xi(A,\mathbf{u}):=\frac{g(\mathbf{u})\kappa_2(A)}{1-g(\mathbf{u})\kappa_2(A)},
  \end{align}
  where $\kappa_2(A)=\|A\|_2\|A^{-1}\|_2$ is the condition number of $A$.
\end{lemma}

\begin{proof}
 Since $A$ is a  symmetric positive definite tridiagonal matrix, it is routine to show by  induction that  the spectral set of $A$ is equal to that of $|A|$. Together with \eqref{eq2}, it yields that
$$\| G\|_2 \le g(\mathbf{u})\| A\|_2.$$
Note that $A$ is nonsingular and $G$ is sufficiently small so that $A+G$ is invertible. Therefore,
\begin{align*}
  \| (I+A^{-1}G)^{-1} \|_2 \le \frac{1}{1-\|A\|_2\|G\|_2}.
\end{align*}
Then it follows from \eqref{fangcheng} that
  \begin{align*}
    \|x-\widehat{x}\|_2&\le \| (I+A^{-1}G)^{-1}\|_2\| A^{-1} \|_2\|G\|_2\|x\|_2\\
    &\le \frac{ \|A^{-1}\|_2\|G\|_2\|x\|_2}{1-\|A^{-1}\|_2\|G\|_2}\le \frac{g(\mathbf{u})\kappa_2(A)}{1-g(\mathbf{u})\kappa_2(A)}\|x\|_2,
  \end{align*}
which completes the proof.
\end{proof}

\begin{theorem}\label{them:jielun}
 Let $A$ be a  symmetric positive definite tridiagonal matrix with
  $\lambda_{\max}(A)$ and $\lambda_{\min}(A)$ being its  largest and  smallest eigenvalue, respectively. Assume that $x$ is the exact solution of the system
    $(A-\alpha I)x=b$
  and $\widehat{x}$ is the computed solution of the system
  $(A-\widetilde{\alpha} I)\widetilde{x}=b$
  by Gaussian elimination, where $\alpha$ is a negative parameter, and
   $\widetilde{\alpha}<0$  is an approximation of $\alpha$ such that $|\widetilde{\alpha}-\alpha|< \varepsilon$ with the  tolerance $\varepsilon$.
  Then, it holds that
 \begin{align}\label{dinglijielun}
   \|x-\widehat{x}\|_2 \le \mathcal{Q}(A,\alpha,\mathbf{u}) \|x\|_2,
 \end{align}
  where
  \begin{align}\label{eq:Q}
  \begin{split}
 \mathcal{Q}(A,\alpha,\mathbf{u}):=&  \frac{g(\mathbf{u})(\lambda_{\max}(A)+|\alpha|+\varepsilon)(\lambda_{\min}(A)+|\alpha|)/(\lambda_{\min}(A)+|\alpha|-\varepsilon)}
   {(\lambda_{\min}(A)+|\alpha|-\varepsilon)-g(\mathbf{u})(\lambda_{\max}(A)+|\alpha|+\varepsilon)}\\
   &+\frac{\varepsilon}{\lambda_{\min}(A)+|\alpha|-\varepsilon}.
   \end{split}
   \end{align}
 Furthermore, if $\widetilde{b}$ is an approximation of $b$ such that
 $\| \widetilde{b}-b \|_2\le \delta$ with  the  tolerance $\delta$, then the computed solution $\widehat{x}$ of the system
 $(A-\widetilde{\alpha} I)\widetilde{x}= \widetilde{b} $  by Gaussian elimination satisfies that
 \begin{align}\label{dingli5}
  \|x-\widehat{x}\|_2\le  
  \mathcal{Q}(A,\alpha,\mathbf{u})\|b\|_2+\frac{ (1+\mathcal{Q}(A,\alpha,\mathbf{u}))\delta }{\lambda_{\min}(A)+|\alpha|}.
 \end{align}

\end{theorem}

\begin{proof}
 We first resort to the conclusion \eqref{dinglijielun}.
 It follows from Lemma \ref{lem:wucha}  that   the computed solution $\widehat{x}$ of the system $(A-\widetilde{\alpha}I)\widetilde{x}=b$ by Gaussian elimination satisfies
 \begin{align*}
   \|\widetilde{x}-\widehat{x}\|_2\le \xi(A-\widetilde{\alpha}I,\mathbf{u})\|\widetilde{x} \|_2.
 \end{align*}
Moreover, it is derived from the perturbation theory of matrices that the exact solution $x$ of  the system $(A-\alpha I)x=b$ satisfies
\begin{align*}
 \| \widetilde{x}-x\|_2\le \varepsilon \| (A-\widetilde{\alpha}I)^{-1} \|_2 \|x\|_2.
\end{align*}
Further, it holds that
\begin{align*}
 \| \widetilde{x}\|_2\le \big(1+\varepsilon\| (A-\widetilde{\alpha}I)^{-1} \|_2 \big)\|x\|_2.
\end{align*}
A combination of  the above estimates leads to
\begin{align}\label{eq3}
   \| x-\widehat{x}\|_2\le \Big[\varepsilon \|(A-\widetilde{\alpha}I)^{-1} \|_2 + \xi(A-\widetilde{\alpha}I,\mathbf{u})\big(1+\varepsilon\| (A-\widetilde{\alpha}I)^{-1} \|_2 \big)\Big]\|x\|_2.
\end{align}
It follows from  the fact that $A$ is a symmetric positive definite matrix and $\widetilde{\alpha}<0$,
$|\widetilde{\alpha}-\alpha|< \varepsilon$ that
\begin{align}
 &\| A- \widetilde{\alpha}I \|_2=\lambda_{\max}(A) + |\widetilde{\alpha}| \le \lambda_{\max}(A) +|\alpha|+\varepsilon,\label{tezhengzhi1}
 \end{align}
 and
 \begin{align}
 & \| (A- \widetilde{\alpha}I)^{-1} \|_2=\frac{1}{\lambda_{\min}(A) + |\widetilde{\alpha}|}\le
 \frac{1}{\lambda_{\min}(A)+|\alpha|-\varepsilon}\label{tezhengzhi2},
\end{align}
which leads to
\begin{align*}
 &\kappa_2(A- \widetilde{\alpha}I)\le \frac{\lambda_{\max}(A)+|\alpha|+\varepsilon}{\lambda_{\min}(A)+|\alpha|-\varepsilon}.
\end{align*}
As a result,
\begin{align}\label{xijielun}
 \xi(A-\widetilde{\alpha}I,\mathbf{u})\le \frac{g(\mathbf{u})(\lambda_{\max}(A)+|\alpha|+\varepsilon )}
 { \lambda_{\min}(A)+|\alpha|-\varepsilon - g(\mathbf{u})(\lambda_{\max}(A)+|\alpha|+\varepsilon)  }.
\end{align}
This  together with \eqref{eq3} leads to  the conclusion \eqref{dinglijielun}.

Then we turn to the conclusion \eqref{dingli5}.
It follows from Lemma \ref{lem:wucha}  that the computed solution $\widehat{x}$ of the system $(A-\alpha I)x^{*}=\widetilde{b}$ by Gaussian elimination satisfies
\begin{align}\label{dinglia}
  \|\widehat{x}-x^{*}\|_2 \le \mathcal{Q}(A,\alpha,\mathbf{u})\|x^{*}\|_2.
\end{align}
A combination of $(A-\alpha I)x^{*}=\widetilde{b}$, \,
 $(A-\alpha I)x=b$ and $\| \widetilde{b}-b \|_2\le \delta$ leads to
\begin{align}\label{dingli3}
  \| x-x^{*}\|_2\le \frac{\delta}{\lambda_{\min}(A)+|\alpha|}, \quad
  \| x^{*}\|_2\le \frac{\delta+ \|b\|_2}{\lambda_{\min}(A)+|\alpha|}.
\end{align}
Together with the triangle inequality, \eqref{dinglia} and \eqref{dingli3}, one obtains that
\begin{align*}
 \|x-\widehat{x}\|_2&\le \mathcal{Q}(A,\alpha,\mathbf{u})\|x^{*}\|_2+
 \frac{\delta}{\lambda_{\min}(A)+|\alpha|}\\
 &\le\mathcal{Q}(A,\alpha,\mathbf{u})\|b\|_2+\frac{(1+\mathcal{Q}(A,\alpha,\mathbf{u}))\delta}{\lambda_{\min}(A)+|\alpha|}.
\end{align*}
The proof is completed.
\end{proof}

\subsection{ Forward error analysis of the feasible ECR algorithm }
In this subsection,  the forward error analysis of the feasible ECR algorithm for the linear  system \eqref{system} will be studied under the following conditions:
\begin{itemize}
 \item[(1)] Both $B$ and $\mathcal{R}_n$ are symmetric positive definite  tridiagonal matrices, and
     all subdiagonal entries of $\mathcal{R}_n$ are non-zero.
 \item[(2)] The eigenvalues of  $B$ and $\mathcal{R}_n$ satisfy
 \begin{align*}
   \lambda_{\text{max}}(\mathcal{R}_n)\le
    \max_{1\le i\le n}|b_i\pm (|a_i|+|a_{i+1}|) |\le 1,  \qquad \lambda_{\text{min}}(\mathcal{R}_n)\geq \mathbf{u}^{1-\epsilon},\\
    \lambda_{\text{max}}(B)\le 1,  \qquad  \lambda_{\text{min}}(B)\geq \mathbf{u}^{1-\epsilon},
 \end{align*}
 where $ 0<\epsilon<1.$ Without loss of generality, we take $\epsilon=\frac{1}{2}$.
\end{itemize}

The above conditions are reasonable  and essential in the following senses:
\begin{itemize}
\item[1,]
 The symmetry requirement of matrix $\mathcal{R}_n$ can be relaxed. Actually,
  for the  non-symmetry tridiagonal matrix $\mathcal{R}_n$ in \eqref{juzhenRn},  if
   $ a_{i+1}c_i >0, \, i=1,2,\cdots,n-1,$ then there exists a diagonal matrix $D=\text{diag}\{d_1,d_2,\cdots,d_n\},$ here
   $$d_1=1, \quad d_j=\sqrt{\frac{a_j}{c_{j-1}}}d_{j-1},\quad j=2,3,\cdots,n,$$
 such that $D^{-1}\mathcal{R}_nD$ is the symmetric tridiagonal matrix and all subdiagonal entries  are non-zero. It follows from the Sturm sequence property \cite{GolubVan2013} that the eigenvalues of any principal submatrix of $\mathcal{R}_n$ are simple, which implies that  the zeros of  matrix polynomial $B_i^{(r)}$ are simple.

\item[2,] By the condition (2) and Lemma \ref{lem:wucha}, it holds that
\begin{align}\label{xi}
\xi(B,\mathbf{u}):=\frac{g(\mathbf{u})\kappa_2(B)}{1-g(\mathbf{u})\kappa_2(B)}<5\sqrt{\mathbf{u}},
\end{align}
which  guarantees that  the round-off error is small enough
 after one step of Gaussian elimination for solving the tridiagonal linear system $Bx=b$.

\item[3,] Let $\mu$ be one zero of  matrix polynomial $B_{i}^{(r)}$,  and $\widetilde{\mu}$  be an approximation of $\mu$ computed by the MOB such that $|\mu-\widetilde{\mu}|<\varepsilon=\frac{15}{2}\mathbf{u}.$
The  ECR algorithm is based on the computation of tridiagonal linear systems of the form $(B-\mu I)x=b$. Assume that $\widehat{x}$ is the computed solution of the system
  $(B-\widetilde{\mu} I)\widetilde{x}=b$ by Gaussian elimination. By Theorem \ref{them:jielun}, one obtains that
  \begin{align*}
    \|x-\widehat{x}\|_2 \le \mathcal{Q}(B,\mu,\mathbf{u}) \|x\|_2,
  \end{align*}
  here $\mathcal{Q}(B,\mu,\mathbf{u})$ is defined in \eqref{eq:Q}.
It  follows from the condition (2) that
    $$2\sqrt{\mathbf{u}}< \lambda_{\min}(B)+|\mu|,$$
   which leads to
     \begin{align*}
       &\varepsilon < \frac{15}{4}\sqrt{\mathbf{u}}(\lambda_{\min}(B)+|\mu|).
     \end{align*}
     This together with $g(\mathbf{u})=\frac{4\mathbf{u}+3\mathbf{u}^2+ \mathbf{u}^3}{1-\mathbf{u}}$ yields that
     \begin{align*}
        &\frac{3.8\varepsilon}{15\sqrt{\mathbf{u}}}< \lambda_{\min}(B)+|\mu|-g(\mathbf{u})(\lambda_{\max}(B)+|\mu|).
     \end{align*}
   Therefore, it holds that
    \begin{align*}
    &\frac{g(\mathbf{u})(\lambda_{\max}(B)\!+\!|\mu|+\varepsilon)}
   {(\lambda_{\min}(B)\!+\!|\mu|\!-\!\varepsilon)\!-\!g(\mathbf{u})(\lambda_{\max}(B)\!+\!|\mu|+\varepsilon)}\\
    &<  \frac{ (1+8\sqrt{\mathbf{u}})g(\mathbf{u})(\lambda_{\max}(B)\!+\!|\mu|)}
   {(\lambda_{\min}(B)+|\mu|)\!-\!g(\mathbf{u})(\lambda_{\max}(B)+|\mu|)},
   \end{align*}
   and
   \begin{align*}
   \frac{\varepsilon}{\lambda_{\min}(B)+|\mu|-\varepsilon}< \frac{3.75}{1-3.75\sqrt{\mathbf{u}}}\frac{\lambda_{\max}(B)+|\mu|}{\lambda_{\min}(B)+|\mu|}.
    \end{align*}
   Finally, it follows from \eqref{xi}, \eqref{dinglijielun} and $\xi(B-\mu I,\mathbf{u})=\frac{g(\mathbf{u})\kappa_2(B-\mu I)}{1-g(\mathbf{u})\kappa_2(B-\mu I)}$ that
 \begin{align}\label{daQ}
  \mathcal{Q}(B,\mu,\mathbf{u})<\frac{1+8\sqrt{\mathbf{u}}+3.75}{1-3.75\sqrt{\mathbf{u}}}\xi(B-\mu I,\mathbf{u})
 <5\,\xi(B-\mu I,\mathbf{u}),
 \end{align}
which guarantees that the round-off error is still small enough after one step of Gaussian elimination for solving $(B-\mu I)x=b$.
\end{itemize}

Let $\mu$ and $\lambda$ be  one zero of matrix polynomials $B_{i}^{(r)}$ and   $B_{i-2^{r-1}}^{(r-1)} B_{i+2^{r-1}}^{(r-1)}$, respectively.
 Let $\widetilde{\mu}<0$ and $\widetilde{\lambda}<0$
  be an approximation of negative parameters $\mu$ and $\lambda$ such that
  $|\widetilde{\mu}-\mu|<\varepsilon$ and $|\widetilde{\lambda}-\lambda|<\varepsilon$ with the tolerance $\varepsilon=\frac{15}{2}\mathbf{u}$, $\widetilde{b}$ be an approximation of $b$ such that
  $\|b-\widetilde{b}\|_2\le \delta<\|b\|_2$ with the tolerance $\delta$.
\begin{theorem}\label{lem:zhongyao}
Assume  that  $x$ is the exact solution of the system $ (B-\mu I)x=(B-\lambda I)b $, and $\widehat{x}$ is the computed solution of the system
$ (B-\widetilde{\mu} I)\widetilde{x}=(B-\widetilde{\lambda} I)b $ by Gaussian elimination.
   Then it holds that
 \begin{align}\label{xiangduiwucha}
   \|\widehat{x}-x\|_2 < \frac{|\mu- \lambda|+|\mu|}{|\mu|}\delta
   + 77\mathbf{u}\frac{\|b\|_2}{|\mu|}.
 \end{align}
\end{theorem}
\begin{proof}

   The system $(B-\mu I)x=(B-\lambda I)b$ is solved by the following two steps:
   \begin{enumerate}[Step 1.]
    \item It is to compute $y$ by the system $(B-\mu I)y=(\mu-\lambda )b$,
    \item It is to compute $x$ by $x=y+b$.
    \end{enumerate}

In  Step 1,  for the computed right-hand side $\widetilde{r}$ effected by the roundoff error of floating-point,  it follows from $|\widetilde{\mu}-\mu|<\varepsilon$, $|\widetilde{\lambda}-\lambda|<\varepsilon$ and
 $\|b-\widetilde{b}\|_2\le \delta<\|b\|_2$ that
\begin{align*}
 \|\widetilde{r}- (\mu-\lambda)b\|_2 \le \mathcal{B}(\mu,\lambda,b,\mathbf{u}),
\end{align*}
where $$\mathcal{B}(\mu,\lambda,b,\mathbf{u})=(|\mu-\lambda|+2\varepsilon) \big[\delta+(\|b\|_2+\delta)\mathbf{u}\big]+2\varepsilon\|b\|_2.$$
By Lemma \ref{lem:wucha}, the computed solution $\widehat{y}$ of the system $(B-\widetilde{\mu}I)\widetilde{y}=\widetilde{r}$ satisfies
\begin{align}\label{eq1}
  \|\widetilde{y}- \widehat{y} \|_2\le \xi(B-\widetilde{\mu}I, \mathbf{u})\| \widetilde{y}\|_2.
\end{align}
 Together with $(B-\mu I)y=(\mu-\lambda)b$ and $|\widetilde{\mu}-\mu|<\varepsilon$, one obtains that
\begin{align}\label{eq2}
  \|\widetilde{y}-y\|_2 \le \|(B-\widetilde{\mu}I)^{-1} \|_2
  \Big( \mathcal{B}(\mu,\lambda,b,\mathbf{u})
  +\varepsilon\|y\|_2 \Big),
\end{align}
which implies
\begin{align}\label{eq3}
 \|\widetilde{y}\|_2\le \big( 1+\varepsilon\| (B-\widetilde{\mu}I)^{-1}\|_2 \big)\|y\|_2
 + \| (B-\widetilde{\mu}I)^{-1} \|_2 \mathcal{B}(\mu,\lambda,b,\mathbf{u}).
\end{align}

In Step 2,  for the computed solution $\widehat{x}$ effected by the roundoff error of floating-point, it holds that
\begin{align}\label{eq4}
\begin{split}
 &\|\widehat{x}-x\|_2
 \le \|\widetilde{y}- \widehat{y}\|_2+ \|\widetilde{y}-y\|_2+ \mathbf{u}\|\widehat{y}\|_2+ \mathbf{u}\|\widetilde{b}\|_2+\delta.
 \end{split}
\end{align}
Furthermore, a combination of \eqref{eq1}, \eqref{eq2} and \eqref{eq3} yields
\begin{align*}
 &\|\widetilde{y}- \widehat{y}\|_2+ \|\widetilde{y}-y\|_2+ \mathbf{u}\|\widehat{y}\|_2\\
 &\le \big[ \mathbf{u} + (1+\mathbf{u})\xi(B-\widetilde{\mu}I,\mathbf{u}) \big]\|\widetilde{y}\|_2+\|(B-\widetilde{\mu}I)^{-1} \|_2  \Big( \mathcal{B}(\mu,\lambda,b,\mathbf{u})
  +\varepsilon\|y\|_2 \Big)\\
  &\le \big[\mathbf{u}+(1+\mathbf{u})\xi(B-\widetilde{\mu}I,\mathbf{u}) \big]
 \big( 1+\varepsilon\|(B-\widetilde{\mu}I)^{-1} \|_2\big)\|y\|_2
 + \varepsilon\|(B-\widetilde{\mu}I)^{-1} \|_2\|y\|_2\\
 &+\big[\mathbf{u}+(1+\mathbf{u})\xi(B-\widetilde{\mu}I,\mathbf{u})+1 \big]
 \|(B-\widetilde{\mu} I)^{-1} \|_2 \mathcal{B}(\mu,\lambda,b,\mathbf{u}).
\end{align*}
In addition, it follows from \eqref{eq:Q}, \eqref{tezhengzhi2} and the fact that  $(B-\mu I)y=(\mu-\lambda )b$ that
\begin{align*}
 &\big[\mathbf{u}+(1+\mathbf{u})\xi(B-\widetilde{\mu}I,\mathbf{u}) \big]\big( 1+\varepsilon\|(B-\widetilde{\mu}I)^{-1} \|_2\big)\|y\|_2+\varepsilon \|(B-\widetilde{\mu}I)^{-1} \|_2\|y\|_2\\
 &\le \big[(1+\mathbf{u})\mathcal{Q}(B,\mu,\mathbf{u})+ \mathbf{u}\big]\|y\|_2\\
  & \le \frac{|\mu-\lambda|}{\lambda_{\min}(B)+|\mu|} \big( \mathbf{u} + (1+\mathbf{u})\mathcal{Q}(B,\mu,\mathbf{u}) \big)\|b\|_2,
\end{align*}
and
\begin{align*}
 &\big[\mathbf{u}(1+\xi(B-\widetilde{\mu}I,\mathbf{u}))+ \xi(B-\widetilde{\mu}I,\mathbf{u})+1 \big]\|(B-\widetilde{\mu}I)^{-1} \|_2 \mathcal{B}(\mu,\lambda,b,\mathbf{u})\\
 &\le (1+\mathbf{u})(1+\xi(B-\widetilde{\mu}I,\mathbf{u}))\|(B-\widetilde{\mu}I)^{-1} \|_2
 \mathcal{B}(\mu,\lambda,b,\mathbf{u})\\
 &\le(1+\mathbf{u})\mathcal{B}(\mu,\lambda,b,\mathbf{u})
 \frac{1+\xi(B-\widetilde{\mu}I,\mathbf{u})}{\lambda_{\min}(B)+|\mu|-\varepsilon}.
\end{align*}

 If the matrix $B$ satisfies conditions (1) and (2), it follows from \eqref{xi} and  \eqref{daQ} that
   \begin{align*}
    &\xi(B-\mu I,\mathbf{u})=\frac{g(\mathbf{u})\kappa_2(B-\mu I)}{1-g(\mathbf{u})\kappa_2(B-\mu I)}< \frac{2g(\mathbf{u})}{|\mu|},\\
      & \mathcal{Q}(B,\mu,\mathbf{u})
         <5\,\xi(B-\mu I,\mathbf{u})<25\sqrt{\mathbf{u}}<0.25,
   \end{align*}
   where $\mathbf{u}$ is sufficiently small such as $\mathbf{u}<10^{-4}$.
   And it gets from  $\varepsilon=\frac{15}{2}\mathbf{u}$ that
    \begin{align*}
     &\lambda_{\min}(B)\!-\!\varepsilon\!-\!g(\mathbf{u})(\lambda_{\max}(B)+|\mu|+\varepsilon)>0,\\
     &\frac{1+\mathbf{u}}{\lambda_{\min}(B)+|\mu|\!-\!\varepsilon\!-\!g(\mathbf{u})(\lambda_{\max}(B)+|\mu|+\varepsilon)}
      <\frac{1}{|\mu|}.
    \end{align*}
As a result,
  \begin{align*}
     &\frac{|\mu-\lambda|}{\lambda_{\min}(B)+|\mu|} \Big[ \mathbf{u} + (1+\mathbf{u})\mathcal{Q}(B,\mu,\mathbf{u}) \Big]\|b\|_2\\
     &< \frac{|\mu-\lambda|}{|\mu|}\Big[ \frac{10g(\mathbf{u})}{|\mu|} + 1.25\mathbf{u} \Big]\|b\|_2
     < \big( \frac{41\mathbf{u}}{|\mu|}+1.25\mathbf{u} \big)\|b\|_2,\\
      \end{align*}
and
  \begin{align*}
    & \frac{(1+\mathbf{u})\mathcal{B}(\mu,\lambda,b,\mathbf{u})}
     {\lambda_{\min}(B)+|\mu|\!-\!\varepsilon\!-\!g(\mathbf{u})
  (\lambda_{\max}(B)+|\mu|+\varepsilon) }\\
  &< \frac{ ( |\mu\!-\!\lambda|+15\mathbf{u})
  (\delta+(\|b\|_2+\delta)\mathbf{u}))+30\mathbf{u}\|b\|_2   }{|\mu|}\\
  &< \frac{|\mu-\lambda|}{|\mu|}\delta + 2\mathbf{u}\|b\|_2+\frac{30\mathbf{u}}{|\mu|}\|b\|_2,
   \end{align*}
 which leads to
    \begin{align*}
  &\|\widehat{x}-x\|_2< \frac{|\mu-\lambda|}{\lambda_{\min}(B)+|\mu|} \big( \mathcal{Q}(B,\mu,\mathbf{u})+\mathbf{u} + \mathbf{u}\mathcal{Q}(B,\mu,\mathbf{u}) \big)\|b\|_2\\
  + & \frac{(1+\mathbf{u})\mathcal{B}(\mu,\lambda,b,\mathbf{u})}
  { \lambda_{\min}(B)+|\mu|-\varepsilon-g(\mathbf{u})(\lambda_{\max}(B)+|\mu|+\varepsilon) }+\mathbf{u}\|b\|_2 + (1+\mathbf{u})\delta\\
   < &\frac{|\mu-\lambda|}{|\mu|}\delta+2\mathbf{u}\|b\|_2
  +\frac{30\mathbf{u}}{|\mu|}\|b\|_2 + (\frac{41\mathbf{u}}{|\mu|}+1.25\mathbf{u})\|b\|_2+ 2\mathbf{u}\|b\|_2+\delta\\
  <& \big[\frac{|\mu-\lambda|}{|\mu|}+1 \big]\delta+\frac{77\mathbf{u}}{|\mu|}\|b\|_2.
 \end{align*}
 The proof is completed.
\end{proof}

 In the reduction phase and the back-substitution phase of the ECR algorithm,
  both implementations  involve solving two  typical   problems of the following form. Take the reduction phase as an example,\newline
 (1)\,for $r=0,1,\ldots,k-2$ and $i=2^{r+1},2\cdot2^{r+1},\ldots, (2^{k-r-1}-1)\cdot2^{r+1}$, compute
     $$ \big( B_i^{(r)}\big)^{-1}B_{i-2^{r-1}}^{(r-1)}B_{i+2^{r-1}}^{(r-1)}p_i^{(r)}.$$
 (2)\,for $r=0,1,\ldots,k-2$ and $i=2^{r+1},2\cdot2^{r+1},\ldots, (2^{k-r-1}-1)\cdot2^{r+1}$,
 compute
 $$\alpha_i^{(r)}\big( B_{i-2^{r-1}}^{(r-1)} \big)^{-1}q_{i-2^r}^{(r)}, \quad \alpha_i^{(r)}=\prod_{j=i-2^{r}+1}^i a_j.$$

The implementation details of these typical   problems have been introduced in Section 3.1.
  The following theorems indicate the implementation process  in this way is stable and the error is controlled.

Let $\{ \mu_i \}_{i=1}^{2^{r+1}-1}$ be the zeros of  matrix polynomial $B_{i}^{(r)} $ such that $0> \mu_1> \mu_2 >\ldots >\mu_{2^{r+1}-1}$, $\{\widetilde{\mu}_i \}_{i=1}^{2^{r+1}-1}$  computed by the MOB  be their corresponding approximations  such that $|\widetilde{\mu}_{i}-\mu_{i}|<\varepsilon$ with the tolerance $\varepsilon=\frac{15}{2}\mathbf{u}$.
Let $\{ \lambda_i \}_{i=1}^{2^{r+1}-2}$ be the zeros of  matrix polynomial $B_{i-2^{r-1}}^{(r-1)}B_{i+2^{r-1}}^{(r-1)}$  such that $0> \lambda_1\geq \lambda_2 \geq \ldots \geq \lambda_{2^{r+1}-2}$,\,$\{\widetilde{\lambda}_i \}_{i=1}^{2^{r+1}-2}$  computed by the MOB
 be their corresponding approximations  such that $|\widetilde{\lambda}_{i}-\lambda_{i}|<\varepsilon$, $\widetilde{b}$ be an approximation of $b$ such that $\|\widetilde{b}-b\|_2< \delta$ with the tolerance $\delta$.
\begin{theorem}\label{dingli1}
 Let $x_1=b,$  compute $x_{2^{r+1}-1}$ by recursively  solving the following
  linear systems
  \begin{align}\label{xitong}
  \begin{split}
    (B-\mu_{\ell+1}I)x_{\ell+1}&= (B-\lambda_{\ell}I)x_{\ell}, \quad  \ell=1,2,\cdots,2^{r+1}-2.
    \end{split}
  \end{align}
Assume  that  the computed solution $\widehat{x}_{2^{r+1}-1}$ holds that
  \begin{align}\label{dingli4}
    \| \widehat{x}_{2^{r+1}-1}-x_{2^{r+1}-1} \|_2\le \frac{|\mu_{2^{r+1}-1}|}{|\mu_1|}
    (\delta+77\mathbf{u}\, \mathcal{C}_1(B_{i}^{(r)}) \|b\|_2),
  \end{align}
  where  $\mathcal{C}_1(B_{i}^{(r)})=\sum_{\ell=1}^{2^{r+1}-2}\frac{1}{|\mu_{\ell+1}|}$.
\end{theorem}

\begin{proof}
  Applying the conclusion of Theorem \ref{lem:zhongyao} to systems \eqref{xitong}  repeatedly, one obtains that
  \begin{align}\label{gongshi}
  \begin{split}
   \| \widehat{x}_{2^{r+1}-1}-x_{2^{r+1}-1} \|_2 <
   &\sum_{\ell=1}^{2^{r+1}-2}\Big[ \frac{77\mathbf{u}}{|\mu_{\ell+1}|}\|x_{\ell}\|_2\cdot\prod_{j=\ell+1}^{2^{r+1}-2}\big( 1+\frac{|\mu_{j+1}-\lambda_{j}|}{|\mu_{j+1}|}   \big)
     \Big]\\
   &+\delta \cdot\prod_{\ell=1}^{2^{r+1}-2}\big( 1+\frac{|\mu_{j+1}-\lambda_{j}|}{|\mu_{j+1}|} \big).
   \end{split}
  \end{align}
 It follows from the fact
  $ (B-\mu_{\ell+1}I)x_{\ell+1}=(B-\lambda_{\ell}I)x_{\ell}=(B-\mu_{\ell+1}I)x_{\ell}
   +(\mu_{\ell+1}-\lambda_{\ell})x_{\ell}$
 that
 \begin{align}\label{gongshi11}
   \|x_{\ell+1}\|_2\le \|x_{\ell}\|_2+\frac{|\mu_{\ell+1}-\lambda_{\ell}|}{\lambda_{\min}(B)+|\mu_{\ell+1}|}\|x_{\ell}\|_2
   \le \prod_{j=1}^{\ell}\big( 1+\frac{|\mu_{j+1}-\lambda_{j}|}{|\mu_{j+1}|} \big)\|x_{1}\|_2.
 \end{align}
In addition, by Lemma \ref{lem:4.2}, it holds that
\begin{align}\label{gongshi22}
  \prod_{\ell=1}^{2^{r+1}-2}\big( 1+\frac{|\mu_{\ell+1}-\lambda_{\ell}|}{|\mu_{\ell+1}|}\big)
  \le \prod_{\ell=1}^{2^{r+1}-2}\frac{|\mu_{\ell+1}|}{|\mu_{\ell}|}=\frac{|u_{2^{r+1}-1}|}{|\mu_1|}.
\end{align}
  Substituting \eqref{gongshi11} and  \eqref{gongshi22} into \eqref{gongshi} leads to the conclusion \eqref{dingli4}.
\end{proof}

\begin{remark}
In the reduction phase of the ECR algorithm, to ensure a stable computing of vector $q_i^{(r)}$, the strategy used is to multiply alternately by the inverse of a factor of $B_i^{(r)}$ and by a factor $B^{(r-1)}_{i-2^{r-1}}B^{(r-1)}_{i+2^{r-1}}$. Since the degree of the corresponding polynomial of matrix  $B^{(r)}_{i}$ is one greater than that of the corresponding polynomial of matrix $B^{(r-1)}_{i-2^{r-1}}B^{(r-1)}_{i+2^{r-1}}$, we leave aside the factor $(B-\mu_1)$ of $B_i^{(r)}$.
 Whether it is the best option or not, it is a feasible choice because of its stability and is  convenient for our theoretical analysis herein.
\end{remark}

Let $\{ \xi_i \}_{i=1}^{2^{r}-1}$ be the zeros of matrix polynomial $B_{i-2^{r-1}}^{(r-1)}$  such that $0> \xi_1\geq \xi_2 \geq \ldots \geq \xi_{2^{r}-1}$,
$\{\widetilde{\xi}_i \}_{i=1}^{2^{r}-1}$  computed by the MOB method
 be their corresponding approximations  such that $|\widetilde{\xi}_{i}-\xi_{i}|<\varepsilon$ with the tolerance $\varepsilon=\frac{15}{2}\mathbf{u}$,\, $\widetilde{b}$ be an approximation of $b$ such that $\|\widetilde{b}-b\|_2< \delta$.
\begin{theorem}\label{dingli2}
 Let $x_0=b$,  compute $x_{2^r-1}$  by recursively
  solving the following linear systems
  \begin{align}\label{houdai}
  \begin{split}
    (B-\xi_{2^r-j}I)x_j &= a_{i-2^r+j}\cdot x_{j-1}, \quad j=1,2,\cdots,2^r-1.
    \end{split}
  \end{align}
  Assume  that  the computed solution $\widehat{x}_{2^{r}-1}$ holds that
   \begin{align}\label{lem:4.5}
     \| \widehat{x}_{2^{r}-1}-x_{2^{r}-1}\|_2 < \mathcal{C}_2(B_{i-2^{r-1}}^{(r-1)}) \big(\delta + 41\mathbf{u}\,\mathcal{C}_3(B_{i-2^{r-1}}^{(r-1)})\|b\|_2\big),
   \end{align}
 where   $\mathcal{C}_2(B_{i-2^{r-1}}^{(r-1)})=\prod_{j=1}^{2^r-1}\frac{|a_{i-j}|}{|\xi_j|}$ and $ \mathcal{C}_3(B_{i-2^{r-1}}^{(r-1)})=\sum_{j=1}^{2^r-1}\frac{1}{|\xi_j|}$.
\end{theorem}

\begin{proof}
Applying the conclusion \eqref{dingli5} of Theorem \ref{them:jielun} to  systems \eqref{houdai} repeatedly, one obtains that
  \begin{align*}
  \|x_{2^r-1}- \widehat{x}_{2^r-1}\|_2
  &\le \sum_{\ell=1}^{2^r-1}\Big[\prod_{j=\ell+1}^{2^r-1}
  (1+\frac{41\mathbf{u}}{|\xi_{2^r-j}|})\frac{|a_{i-2^r+j}|}{\lambda_{\min}(B)+|\xi_{2^r-j}|}\Big]
  \frac{41\mathbf{u}}{|\xi_{2^r-\ell}|}\|x_{\ell}\|_2\\
  &+ \delta \prod_{j=1}^{2^r-1}
  (1+\frac{41\mathbf{u}}{|\xi_{2^r-j}|})\frac{|a_{i-2^r+j}|}{\lambda_{\min}(B)+|\xi_{2^r-j}|}.
  \end{align*}
  Moreover, for $1\le j \le 2^r-1$, it holds that
  \begin{align*}
  \frac{1+\frac{41\mathbf{u}}{|\xi_{2^r-j}|}}{\lambda_{\min}(B)+|\xi_{2^r-j}|}
 =\frac{1}{\xi_{2^r-j}}\frac{41\mathbf{u}+|\xi_{2^r-j}|}{\lambda_{\min}(B)+|\xi_{2^r-j}|}
 <\frac{1}{|\xi_{2^r-j}|},
 \end{align*}
 which implies
 \begin{align}\label{gongshi5}
    \|x_{2^r-1}- \widehat{x}_{2^r-1}\|_2\le
    \sum_{\ell=1}^{2^r-1}\Big(\prod_{j=\ell+1}^{2^r-1}
  \frac{|a_{i-2^r+j}|}{|\xi_{2^r-j}|}\Big)
  \frac{41\mathbf{u}}{|\xi_{2^r-\ell}|}\|x_{\ell}\|_2+ \delta \prod_{j=1}^{2^r-1}
 \frac{|a_{i-2^r+j}|}{|\xi_{2^r-j}|}.
 \end{align}
 It follows from $(B-\xi_{2^r-\ell}I)x_\ell = a_{i-2^r+\ell}\cdot x_{\ell-1}$ that
 \begin{align}\label{gongshi4}
   \|x_{\ell}\|_2\le \frac{|a_{i-2^r+\ell}|}{\lambda_{\min}(B)+|\xi_{2^r-\ell}|}\|x_{\ell-1}\|_2
   \le \big( \prod_{j=1}^{\ell}\frac{|a_{i-2^r+j}|}{|\xi_{2^r-j}|} \big)\|b\|_2.
 \end{align}
Finally, substituting \eqref{gongshi4} into \eqref{gongshi5} leads to the conclusion \eqref{lem:4.5}.
\end{proof}

\begin{remark}
 To make a comparison, the  backward error of the column principal  element Gaussian elimination method for the linear system with the matrix $B_{i}^{(r)}$ is investigated. As is shown in Lemma
 \ref{forwarderror} below,
 the upper bound of the backward error is much larger than that of the forward error of  recursively solving the tridiagonal linear systems with the factor of $B_{i}^{(r)}$ \eqref{xitong} in Theorem \ref{dingli1} and \eqref{houdai} in Theorem \ref{dingli2}.
 \begin{lemma}[P67, \cite{Xu2013}]\label{forwarderror}
   Let $A\in \RR^{n\times n}$ be a nonsingular matrix and $1.01n\mathbf{u}\le 0.01$. The computed solution $\widehat{x}$ of the linear system $Ax=b$ by the column principal  element Gaussian elimination method  satisfies that
   $$(A+\delta A)\widehat{x}=b.$$
Then it holds that
 $$\frac{\|\delta A\|_{\infty}}{\|A\|_{\infty}}\le 4.09n^3\rho \mathbf{u}.$$
It has been proved that $\rho\le 2^{n-1}$ in theory and the upper bound $2^{n-1}$ can be reached.
 \end{lemma}
\end{remark}

\subsection{ The ECR algorithm for special systems }

 It is to use the ECR algorithm for  linear systems resulting from the Legendre spectral-Galerkin method  for the Poisson equation on a square domain.

Let $b_j$ and $a_j$ be diagonal entries and subdiagonal entries ($j=1,2,\cdots,n$, $a_1=0$ and $a_{n+1}=0$) of the following symmetric tridiagonal matrix  $M_i (i=1,2)$ \cite{Shen1995} of order $n$, respectively,
\begin{enumerate}[(1)]
\item
for matrix $ M_1,\, b_{i}=\frac{2}{(4i-3)(4i+1)}, \quad a_{i+1}=\frac{-1}{\sqrt{(4i+3)(4i-1)}(4i+1)},$ \qquad $i=1,2,\cdots,n,$
\item
for matrix $M_2,\, b_{i}=\frac{2}{(4i-1)(4i+3)}, \quad  a_{i+1}=\frac{-1}{\sqrt{(4i+5)(4i+1)}(4i+3)},$ \qquad $i=1,2,\cdots,n.$
\end{enumerate}

It is obvious that conditions (1) and (2) in Section 4.3 hold for  $M_1$ and $M_2$. In what
follows, $\mathcal{C}_1$ in Theorem \ref{dingli1} and $\mathcal{C}_2, \mathcal{C}_3$ in Theorem \ref{dingli2} are investigated. For this purpose, the
DETGTRI algorithm \cite{ElMikkawy2004} is presented which is an efficient computational method for evaluating the  determinant of a tridiagonal matrix of order $n$ with only cost $\mathcal{O}(n)$ .

To compute the  determinant of a tridiagonal matrix of order $n$ such as $\mathcal{R}_n$ shown in \eqref{juzhenRn}, it is proceeded as follows:
\begin{enumerate}[Step 1.]
\item Introduce an  additional $n$-dimensional vector $\mathbf{g}=(g_1,g_2,\cdots,g_n)$ by
\begin{align}\label{Detgtri}
 g_i = \begin{cases}
  b_1, & \text{if} \,\,i=1,\\
  b_i-\frac{a_ic_{i-1}}{g_{i-1}}, & \text{if} \,\,i=2,3,\cdots,n.
 \end{cases}
 \end{align}
 Use \eqref{Detgtri} to compute  the $n$ components of the vector $\mathbf{g}$. If $g_i=0$ for any $i\le n$, set $g_i=x$ ($x$ is just a symbolic name) and continue to compute $g_{i+1}, g_{i+2},\cdots, g_n$ in terms of $x$ by using \eqref{Detgtri}.
 \item The  product $\prod_{i=1}^ng_i$\,(this product is a polynomial in $x$) evaluated at $x=0$ is equal to the determinant of the tridiagonal matrix $\mathcal{R}_n$.
\end{enumerate}

Note that  the  product $\prod_{i=1}^ng_i$ is the exact determinant of a tridiagonal matrix in theory.

\begin{theorem}\label{juzhenhanglieshi}
   Given any principal submatrix of the symmetric  tridiagonal matrix $M_i (i=1,2)$ as follows
    \begin{align*}
  R_{r-t+1}=
\begin{bmatrix}
 b_{t} & a_{t+1} &     &      \\
 a_{t+1} & b_{t+1} & a_{t+2} &      \\
          & \ddots &\ddots &\ddots \\
     &      &  a_{r-1} & b_{r-1} & a_{r} \\
     &      &         & a_{r} & b_{r}
\end{bmatrix},
\quad r-t+1\geq 3.
\end{align*}
 Let $a_t=b_t$. If $\frac{b_{k+1}}{|a_{k+1}|}-\frac{|a_{k+1}|}{|a_{k}|}>1$ for $k=t,t+1,\cdots,r-1$,  it holds that
\begin{align*}
 b_t \prod_{i=t+1}^r |a_i| < \det(R_{r-t+1}).
\end{align*}
\end{theorem}

\begin{proof}
  According to the DETGTRI algorithm, the  vector $\mathbf{g}=(g_t,g_{t+1},\cdots,g_{r})$ reads
  \begin{align*}
 g_i = \begin{cases}
  b_t, & \text{if} \,\,i=t,\\
  b_i-\frac{a_i a_{i}}{g_{i-1}}, & \text{if} \,\,i=t+1,t+2,\cdots,r.
 \end{cases}
 \end{align*}
 We will prove $g_{k}\geq |a_k|$ for $k=t,t+1,\cdots,r$ by induction on $k$ under the following condition
 \begin{align}\label{tiaojian}
 \frac{b_{k+1}}{|a_{k+1}|}-\frac{|a_{k+1}|}{|a_{k}|}>1.
 \end{align}
It is obvious that  $g_k=a_k$ for $k=t$. For $k=t+1,$ it follows from \eqref{tiaojian} that
  \begin{align*}
    \frac{g_{t+1}}{|a_{t+1}|}=\frac{b_{t+1}}{|a_{t+1}|}-\frac{|a_{t+1}|}{g_t}>1,
  \end{align*}
  which leads to  $g_{t+1}> |a_{t+1}|.$ Assume $g_k \geq |a_{k}|$, together with \eqref{tiaojian},  it holds for
  $k+1$ that
  \begin{align*}
    \frac{g_{k+1}}{|a_{k+1}|}=\frac{b_{k+1}}{|a_{k+1}|}-\frac{|a_{k+1}|}{g_k}\geq
    \frac{b_{k+1}}{|a_{k+1}|}-\frac{|a_{k+1}|}{|a_k|}>1,
  \end{align*}
  which proves the conclusion $g_{k}\geq |a_k|$ for $k=t,t+1,\cdots,r$. As a result,
  \begin{align*}
     \det(R_{r-t+1})= \prod_{i=t}^r g_i > b_t \prod_{i=t+1}^r |a_i|.
  \end{align*}
  The proof is completed.
\end{proof}

For  matrices $M_1$ and $M_2$,  $\mathcal{C}_1(B_{i}^{(r)})$ in Theorem \ref{dingli1},  $\mathcal{C}_2(B_{i-2^{r-1}}^{(r-1)})$ and $\mathcal{C}_3(B_{i-2^{r-1}}^{(r-1)})$ in Theorem \ref{dingli2} are evaluated as follows:
\begin{enumerate}[(1)]
\item For $M_1$, $r=1,2,\cdots,k-1$, $i=j\cdot2^r,$ $j=1,2,\cdots,2^{k-r-1},$
 $$\mathcal{C}_1(B_{i}^{(r)})< \frac{2^{r+1}-2}{|\lambda_{\min}(-\mathcal{R}_i^{(r)})|}, \quad  \qquad
 \mathcal{C}_3(B_{i-2^{r-1}}^{(r-1)})< \frac{2^r-1}{|\lambda_{\min}(-\mathcal{R}_{i-2^{r-1}}^{(r-1)})|}.$$
It is easy to check that for any principal submatrix $R_{\ell-t+1}$ of $M_1$, it holds that
    \begin{align*}
       \det(R_{\ell-t+1})>
       \begin{cases}
         b_t\cdot\prod\limits_{i=t+1}^{\ell} |a_i|/1.1, \quad &\text{if}\,\, t=1,\ell=3,\\
         b_t\cdot\prod\limits_{i=t+1}^{\ell} |a_i|,  \quad &\text{else},
       \end{cases}
    \end{align*}
    which leads to
            $$\mathcal{C}_2(B_{i-2^{r-1}}^{(r-1)})< \frac{1.1}{|b_{i-2^r+1}|},\quad b_{i-2^r+1}=\frac{2}{(4i-2^{r+2}+1)(4i-2^{r+2}+5)},$$
            where $b_{i-2^r+1}$ is the element in the first row and the first column of matrix $-\mathcal{R}_{i-2^{r-1}}^{(r-1)}$.
 \item For $M_2$, $r=1,2,\cdots,k-1$, $i=j\cdot2^r,$ $j=1,2,\cdots,2^{k-r-1},$
 $$\mathcal{C}_1(B_{i}^{(r)})< \frac{2^{r+1}-2}{|\lambda_{\min}(-\mathcal{R}_i^{(r)})|}, \quad  \qquad
 \mathcal{C}_3(B_{i-2^{r-1}}^{(r-1)})< \frac{2^r-1}{|\lambda_{\min}(-\mathcal{R}_{i-2^{r-1}}^{(r-1)})|}.$$
 It is easy to verify that
       $\frac{b_{k+1}}{|a_{k+1}|}-\frac{|a_{k+1}|}{|a_{k}|}>1$ for $k=1,2,\cdots,n$.
       Then it follows from Theorem \ref{juzhenhanglieshi} that
       for any principal submatrix $R_{\ell-t+1}$ of $M_2$,
       \begin{align*}
       \det(R_{\ell-t+1})> b_t\cdot\prod\limits_{i=t+1}^{\ell} |a_i|,
    \end{align*}
    which leads to
    $$\mathcal{C}_2(B_{i-2^{r-1}}^{(r-1)})< \frac{1}{|b_{i-2^r+1}|},\quad b_{i-2^r+1}=\frac{2}{(4i-2^{r+2}+3)(4i-2^{r+2}+7)}, $$
 where $b_{i-2^r+1}$ is the element in the first row and the first column of matrix $-\mathcal{R}_{i-2^{r-1}}^{(r-1)}$.
\end{enumerate}

\begin{appendix}
\section{The proof of \eqref{gongshi1} and \eqref{gongshi2} }
We first  prove \eqref{gongshi1}. Let
\begin{align*}
& \mathbf{\widehat{a}}_{i-2^{t-1}} = (0,\cdots,0,a_{i-2^{t-1}})\in \RR^{2^{t-1}-1},\\
&\mathbf{\widehat{c}}_{i-2^{t-1}} =  (c_{i-2^{t-1}},0,\cdots,0)\in \RR^{2^{t-1}-1},\\
&\mathbf{\widehat{c}}_{i-2^{t-1}-1} =  (0,\cdots,0,c_{i-2^{t-1}-1})^{\mathrm{ T }}\in\RR^{2^{t-1}-1},\\
&\mathbf{\widehat{a}}_{i-2^{t-1}+1} = (a_{i-2^{t-1}+1},0,\cdots,0)^{\mathrm{ T }}\in \RR^{2^{t-1}-1},
\end{align*}
and
\begin{align*}
  &\mathcal{A}_{1}=
\begin{bmatrix}
xI-\mathcal{R}_{i-3\cdot 2^{t-2}}^{(t-2)} & \mathbf{\widehat{c}}_{i-2^{t-1}-1}\\
\mathbf{\widehat{a}}_{i-2^{t-1}} &  x+b_{i-2^{t-1}}
\end{bmatrix},
\mathcal{A}_{2}=
\begin{bmatrix}
0 & \mathbf{\widehat{a}}_{i-2^{t-1}+1} \\
0 & 0
\end{bmatrix},
  \mathcal{A}_{3}=
  \begin{bmatrix}
    0 & 0 \\
 \mathbf{\widehat{c}}_{i-2^{t-1}} & 0
   \end{bmatrix},\quad\\
 &\mathcal{A}_{4}=
 \begin{bmatrix}
    xI-\mathcal{R}_{i- 2^{t-2}}^{(t-2)} & 0 \\
   0 & xI-\widehat{\mathcal{R}}_{i-2^{t-2}}^{(t-2)}
 \end{bmatrix},\quad
 \mathcal{A}_{5}=
\begin{bmatrix}
    xI-\widehat{\mathcal{R}}_{i-2^{t-2}}^{(t-2)} & 0 \\
   0 & xI-\mathcal{R}_{i- 2^{t-2}}^{(t-2)}\\
\end{bmatrix},\\
&\mathcal{A}_{22}=
\begin{bmatrix}
0 & 0 \\
0 & \mathbf{\widehat{a}}_{i-2^{t-1}+1}
\end{bmatrix},\quad
  \mathcal{A}_{33}=
  \begin{bmatrix}
    0 & 0 \\
    0&\mathbf{\widehat{c}}_{i-2^{t-1}}
   \end{bmatrix}.
\end{align*}
By interchanging rows and columns of a matrix, it yields
\begin{align*}
& f_{i-2^{t-1}}^{(t-1)}(x)
\det\Big( xI - \widehat{\mathcal{R}}_{i-2^{t-2}}^{(t-2)}\Big)
- f_{i-2^{t-2}}^{(t-2)}(x)\det\Big(xI - \widehat{\mathcal{R}}_{i-2^{t-1}}^{(t-1)}\Big)\\
=&\det\Big(
  \begin{bmatrix}
        \mathcal{A}_{1} & \mathcal{A}_{3} \\
           \mathcal{A}_{2} & \mathcal{A}_{4} \\
\end{bmatrix}\Big)
-
\det\Big(
  \begin{bmatrix}
        \mathcal{A}_{1} & \mathcal{A}_{3} \\
         \mathcal{A}_{2} & \mathcal{A}_{5} \\
\end{bmatrix}\Big)
= \det\Big(
  \begin{bmatrix}
        \mathcal{A}_{1} & \mathcal{A}_{3}  \\
        \mathcal{A}_{2}  & \mathcal{A}_{4} \\
\end{bmatrix}\Big)
-
\det\Big(
  \begin{bmatrix}
        \mathcal{A}_{1} & \mathcal{A}_{33} \\
        \mathcal{A}_{22} & \mathcal{A}_{4} \\
\end{bmatrix}\Big),
\end{align*}
Furthermore, one obtains that
\begin{align*}
& \det\Big(
  \begin{bmatrix}
        \mathcal{A}_{1} & \mathcal{A}_{3}  \\
        \mathcal{A}_{2}  & \mathcal{A}_{4} \\
\end{bmatrix}\Big)
-
\det\Big(
  \begin{bmatrix}
        \mathcal{A}_{1} & \mathcal{A}_{33} \\
        \mathcal{A}_{22} & \mathcal{A}_{4} \\
\end{bmatrix}\Big)\\
=& \det\Big(
  \begin{bmatrix}
        \mathcal{A}_{1} & \mathcal{A}_{3}  \\
        \mathcal{A}_{2}  & \mathcal{A}_{4} \\
\end{bmatrix}\Big)
-
\det\Big(
  \begin{bmatrix}
        \mathcal{A}_{1} & \mathcal{A}_{3}  \\
        \mathcal{A}_{22}  & \mathcal{A}_{4} \\
\end{bmatrix}\Big)
+
\det\Big(
  \begin{bmatrix}
        \mathcal{A}_{1} & \mathcal{A}_{3}  \\
        \mathcal{A}_{22}  & \mathcal{A}_{4} \\
\end{bmatrix}\Big)
-
\det\Big(
  \begin{bmatrix}
        \mathcal{A}_{1} & \mathcal{A}_{33} \\
        \mathcal{A}_{22} & \mathcal{A}_{4} \\
\end{bmatrix}\Big),\\
=& \det\Big(\begin{bmatrix}
xI- \mathcal{R}_{i-3\cdot 2^{t-2}}^{(t-2)} & 0 &0 & 0 \\
\mathbf{\widehat{a}}_{i-2^{t-1}} &  0 &  \mathbf{\widehat{c}}_{i-2^{t-1}} & 0 \\
0 & \mathbf{\widehat{a}}_{i-2^{t-1}+1} &  xI- R_{i- 2^{t-2}}^{(t-2)} & 0 \\
0 & -\mathbf{\widehat{a}}_{i-2^{t-1}+1} & 0 & xI- \widehat{\mathcal{R}}_{i-2^{t-2}}^{(t-2)}\\
 \end{bmatrix} \Big)\\
 -&\det\Big(\begin{bmatrix}
xI-  \mathcal{R}_{i-3\cdot 2^{t-2}}^{(t-2)} & \mathbf{\widehat{c}}_{i-2^{t-1}-1} &0 & 0 \\
0 &  0 &  \mathbf{\widehat{c}}_{i-2^{t-1}} & -\mathbf{\widehat{c}}_{i-2^{t-1}} \\
0 & 0 &  xI- \mathcal{R}_{i- 2^{t-2}}^{(t-2)} & 0 \\
0 & \mathbf{\widehat{a}}_{i-2^{t-1}+1} & 0 & xI- \widehat{\mathcal{R}}_{i-2^{t-2}}^{(t-2)}\\
 \end{bmatrix} \Big)\\
=&a_{i-2^{t-1}+1} c_{i-2^{t-1}}f_{i-3\cdot 2^{t-2}}^{(t-2)}(x)
 \det\Big(xI- \mathcal{R}_{i- 2^{t-2}}^{(t-2)}\Big)
  \det\Big(xI- \breve{\mathcal{R}}_{i- 2^{t-2}}^{(t-2)}\Big)\\
  -&a_{i-2^{t-1}+1} c_{i-2^{t-1}}f_{i-3\cdot 2^{t-2}}^{(t-2)}(x)
 \det\Big(xI- \widehat{\mathcal{R}}_{i-2^{t-2}}^{(t-2)}\Big)
  \det\Big(xI-\widetilde{\mathcal{R}}_{i- 2^{t-2}}^{(t-2)}\Big),
\end{align*}

where
\begin{align*}
& \breve{R}_{i-2^{t-2}}^{(t-2)}= -\mathcal{R}_{n}\big([i-2^{t-1}:i-2;i-2^{t-1}:i-2]\big).
\end{align*}
Then it follows from Lemma \ref{lemma} that
\begin{align*}
 &\det\left(xI- \mathcal{R}_{i- 2^{t-2}}^{(t-2)}\right)
  \det\left(xI-   \breve{\mathcal{R}}_{i- 2^{t-2}}^{(t-2)}\right) -
   \det\left(xI- \widehat{\mathcal{R}}_{i-2^{t-2}}^{(t-2)}\right)
  \det\left(xI-  \widetilde{\mathcal{R}}_{i- 2^{t-2}}^{(t-2)}\right)\\
   &= -\prod\limits_{j=i-2^{t-1}+1}^{i-2} a_{j+1}c_j.
\end{align*}
As a result,
\begin{align*}
&\det\left(xI - \mathcal{R}_{i-2^{t-1}}^{(t-1)}\right)
\det\left( xI - \widehat{\mathcal{R}}_{i-2^{t-2}}^{(t-2)} \right)
- \det\left(xI - \widehat{R}_{i-2^{t-1}}^{(t-1)}\right)
\det\left( xI -  \mathcal{R}_{i-2^{t-2}}^{(t-2)}\right)\\
&= -\det\left(xI- \mathcal{R}_{i-3\cdot 2^{t-2}}^{(t-2)}\right)
\prod\limits_{j=i-2^{t-1}}^{i-2} a_{j+1}c_j,
\end{align*}
which proves \eqref{gongshi1}.
Similar process leads to the conclusion below
\begin{align*}
&\det\left(xI - \mathcal{R}_{i+2^{t-1}}^{(t-1)}\right)
\det\left( xI -  \widetilde{\mathcal{R}}_{i+2^{t-2}}^{(t-2)}\right)
- \det\left(xI - \widetilde{\mathcal{R}}_{i+2^{t-1}}^{(t-1)}\right)
\det\left( xI -  \mathcal{R}_{i+2^{t-2}}^{(t-2)}\right)\\
&= -\det\left(xI- \mathcal{R}_{i+3\cdot 2^{t-2}}^{(t-2)}\right)
\prod\limits^{i+2^{t-1}-1}_{j=i+1} a_{j+1}c_j,
\end{align*}
which proves \eqref{gongshi2}.

\end{appendix}

\section*{Acknowledgements}
The work of the third author is supported in part by the National Natural Science Foundation of China (No. 12101325).


\begin{thebibliography}{4}

\bibitem{BiniMeini2009}
 D. A. Bini, B. Meini, The cyclic reduction algorithm: from Poisson equation to stochastic processes and beyond,  Numer. Algorithms, 51 (2009), pp. 23-60.
 
 \bibitem{Bank1975}
R. E. Bank, Marching algorithms for elliptic boundary value problems, Ph.D. thesis, Harvard
Univ., Cambridge, Mass., 1975.

\bibitem{BankRose1975}
R. E. Bank and D. J. Rose, An $\mathcal{O}(n^2)$ method for solving constant coefficient boundary value problems in two dimensions, SIAM J. Numer. Anal., 12 (1975), pp. 529-540.

\bibitem{BuzbeeGolub1970}
B. L. Buzbee, G. H. Golub, and C. W. Nielson, On direct methods for solving Poisson's
equations, SIAM J. Numer. Anal., 7 (1970), pp. 627-656.

\bibitem{BuzbeeDorr1974}
B. L. Buzbee and F. W. Dorr, The direct solution of the biharmonic equation on rectangular
regiots and the Poisson equation on irregular regions, Ibid., 11 (1974), pp. 753-763.

\bibitem{BuzbeeDorr1971}
B. L. Buzbee, F. W. Dorr,  J. A. George and G. H. Golub, The direct solutions of the
discrete Poisson equation on irregular regions, SIAM  J. Numer. Anal., 8 (1971), pp. 722-736.


\bibitem{Barth1967}
W. Barth,  R. S. Martin, and J. H. Wilkinson. Calculation of the eigenvalues of a symmetric tridiagonal matrix by the method of bisection,  Numer. Math., 9 (1967), pp. 386-393.



\bibitem{Legendre1}
X. H. Diao, J. Hu, and S. N. Ma, Preconditioned Legendre spectral Galerkin methods for the non-separable elliptic equation, 91 (2022), pp. 1-27.

\bibitem{ElMikkawy2004}
M. E. A.  El-Mikkawy. A fast algorithm for evaluating nth order tri-diagonal determinants, J. Comput. Appl. Math., 166 (2004), pp. 581-584.

\bibitem{Higham1990}
N. J. Higham, Bounding the Error in Gaussian Eimination for Tridiagonal Systems, SIAM J. Matrix Anal. A., 11 (1990), pp. 521-530.

\bibitem{HaidvogelZang1979}
D. B. Haidvogel and T. A. Zang, The accurate solution of Poisson's equation by expansion
in Chebyshev polynomials , J. Comput. Phys., 30 (1979), pp. 167-180.

\bibitem{GolubVan2013}
 G. H. Golub and C.F. Van Loan, Matrix Computations, The Johns Hopkins University Press,
Baltimore, MD, 4nd ed., 2013.

\bibitem{Stewart1990}
  D. P. O'leary and  G. W. Stewart, Computing the eigenvalues and eigenvectors of symmetric arrowhead matrices, J. Comput. Phys., 90 (1990), pp. 497-505.
  

\bibitem{Swarztrauber1974JCP}
P. N. Swarztrauber, The direct solution of the discrete Poisson equation on the surface of a
sphere, J. Comput. Phys., 15 (1974), pp. 46-54.

\bibitem{Swarztrauber1973}
 P. N. Swarztrauber and R. A. Sweet. The direct solution of the discrete Poisson equation on
a disk, Ibido, 10 (1973), pp. 900-907.


\bibitem{Swarztrauber1974}
 P. N. Swarztrauber, A direct method for the discrete solution of separable elliptic equations, SIAM J. Numer. Anal., 11 (1974), pp. 1136-1150.

 \bibitem{Shen1995}
 J. Shen, On fast direct poisson solver, inf-sup constant and iterative Stokes solver by Legendre Galerkin method, J. Comput. Phys.,  116 (1995), pp. 184-188.

 \bibitem{Shen1994Legendre}
  J. Shen, Efficient spectral-Galerkin method I. Direct solvers for second- and fourth-order equations using Legendre polynomials, SIAM J. Sci. Comput., 15 (1994), pp. 1489-1505.


\bibitem{Swarztrauber1975}
P. N. Swarztrauber and R. A.  Sweet, Efficient FORTRAN Subprograms for the Solution of Elliptic Equations,
NCAR Technical Report TN/IA-109, National Center for Atmospheric Research, 1975.

\bibitem{saad2003}
  Y. Saad, Iterative methods for sparse linear systems, Society for Industrial and Applied Mathematics,
  Philadelphia, 2nd ed., 2003.


 \bibitem{Swarztrauber1977}
   P. N. Swarztrauber, The methods of Cyclic reduction, Fourier analysis and the FACR algorithm for the discrete solution of Poisson's equation on a rectangle. SIAM Rev., 19 (1977), pp. 490-501.


\bibitem{WalterGolub1997}
  G. Walter and G. H. Golub, Cyclic Reduction-History and Applications, scientific computing, 1997.

 \bibitem{Xu2013}
  S. F. Xu, L. Gao, and P. W Zhang, Numerical Linear Algebra,  Peking University Press, Beijing,  2nd ed., 2013.

%
%
%
%
%
%
%
%
%
%
%
%
%


\end{thebibliography}
\end{document}